\documentclass[12pt]{amsart}
 \usepackage[dvips]{graphicx}

\usepackage{amsmath,graphics}
\usepackage{amsfonts,amssymb}

\theoremstyle{plain}
\newtheorem*{theorem*}{Theorem}
\newtheorem*{lemma*} {Lemma}
\newtheorem*{corollary*} {Corollary}
\newtheorem*{proposition*} {Proposition}
\newtheorem*{conjecture*}{Conjecture}

\newtheorem{theorem}{Theorem}[section]
\newtheorem{lemma}[theorem]{Lemma}
\newtheorem{corollary}[theorem]{Corollary}
\newtheorem{proposition}[theorem]{Proposition}

\theoremstyle{definition}

\newcommand{\ord}{\operatorname{ord}}

\newcommand{\ann}{\operatorname{Ann}}

\newcommand{\id}{\operatorname{id}}

\theoremstyle{remark}

\theoremstyle{definition}
\newtheorem{defn}[theorem]{Definition}

\def\s{\sigma}

\def\gl{\mbox{GL}}
\def\Q{\Bbb{Q}}
\def\F{\Bbb{F}}

\def\Z{\Bbb{Z}}
\def\R{\Bbb{R}}

\def\N{\Bbb{N}}
\def\genus{\mbox{genus}}

\def\part{\partial}

\def\a{\alpha}

\def\tor{\mbox{Tor}}
\def\bp{\begin{pmatrix}}

\def\sm{\setminus}
\def\ep{\end{pmatrix}}
\def\bn{\begin{enumerate}}

\def\Hom{\mbox{Hom}}

\def\rank{\mbox{rank}}

\def\en{\end{enumerate}}
\def\ba{\begin{array}}
\def\ea{\end{array}}

\def\L{\Lambda}

\def\at{\overline{\a}}
\def\s{\sigma}
\def\a{\alpha}

\def\ti{\tilde}

\def\vph{\varphi}
\def\fr12{\frac{1}{2}}

\def\im{\mbox{Im}}

\def\ker{\mbox{Ker}}

  \def\hom{\mbox{Hom}} \def\Ext{\mbox{Ext}}

\def\ext{\mbox{Ext}}       \def\dim{\mbox{dim}} 
 \def\xpm {[x_1^{\pm 1}]}
 \def\xipm {[x_i^{\pm 1}]}
 \def\ypm {[x_2^{\pm 1}]}
 \def\xypm {[x_1^{\pm 1},x_2^{\pm 1}]}
 \def\be{\begin{equation} }
 \def\ee{\end{equation}}

\def\twialex{\Delta_{\phi}^{\a}(t)}
\def\twialexi{\Delta_{\phi}^{\a,i}(t)}
\def\twialexzero{\Delta_{\phi}^{\a,0}(t)}
\def\twialextwo{\Delta_{\phi}^{\a,2}(t)}

 \def\tpm{[t^{\pm 1}]}

\def\ft{\F[t^{\pm 1}]}
\def\f{\F}
\def\GL{\mbox{GL}}
\def\glfk{\GL(\F,k)}
\def\fk{\F^k}

\def\tnphi{||\phi||_T}

\def\degtwialex{\deg\left(\Delta_{\phi}^{\a}(t)\right)}
\def\degtwialextwo{\deg\left(\Delta_{\phi}^{\a,2}(t)\right)}
\def\degtwialexzero{\deg\left(\Delta_{\phi}^{\a,0}(t)\right)}

\begin{document}

\title{Twisted Alexander norms give lower bounds on the Thurston norm}
\author{Stefan Friedl and Taehee Kim}
\date{\today} \address{Rice University, Houston, Texas, 77005-1892}
\email{friedl@math.rice.edu}\address{ Department of Mathematics,
Konkuk University, Hwayang-dong, Gwangjin-gu, Seoul 143-701, Korea}
\email{tkim@konkuk.ac.kr}
\def\subjclassname{\textup{2000} Mathematics Subject Classification}
\expandafter\let\csname
subjclassname@1991\endcsname=\subjclassname
\expandafter\let\csname
subjclassname@2000\endcsname=\subjclassname \subjclass{Primary
57M27; Secondary 57N10} \keywords{Thurston norm, Twisted Alexander
norm, 3-manifolds}

\begin{abstract}
We introduce twisted Alexander norms of a compact connected orientable 3-manifold with first
Betti number bigger than one, generalizing norms of McMullen and Turaev. We show that twisted
Alexander norms give lower bounds on the Thurston norm of a 3-manifold. Using these we
completely determine the Thurston norm of many 3-manifolds which are not  determined by norms of
McMullen and Turaev.
\end{abstract}

\maketitle

%
%
%
%
%

%

\section{Introduction}
\label{sec:introduction}
Let $M$ be a 3-manifold. Throughout the paper we will assume that all 3-manifolds are compact,
connected and orientable. Let $\phi\in H^1(M;\Z)$. There exists a (possibly disconnected) properly
embedded surface $S$ which represents a homology class which is dual to $\phi$. (We also say that
$S$ is \emph{dual} to $\phi$.) The \emph{Thurston norm} of $\phi$ is now defined as
 \[
||\phi||_{T,M}=\min \{ -\chi(\hat{S})\, | \, S \subset M \mbox{ properly embedded surface dual to
}\phi\}
\] where $\hat{S}$ denotes the result of discarding all connected components of $S$ with positive Euler
characteristic.
If the manifold $M$ is clear, we will just write $||\phi||_T$.

Thurston \cite{Th86} introduced $||-||_T$ in a preprint in 1976.
He proved that the Thurston norm
on $H^1(M;\Z)$ is homogeneous and convex (that is, for
$\phi,\phi_1,\phi_2\in H^1(M;\Z)$ and $k\in \N$, $||k\phi||_T=k
||\phi||_T$ and $||\phi_1+\phi_2||_T\leq
||\phi_1||_T+||\phi_2||_T$). He also showed that the Thurston norm
can be extended to a seminorm on $H^1(M;\R)$ and that the Thurston norm ball
(which is the set of $\phi\in H^1(M;\R)$ with $||\phi||_T \le 1$)
is a (possibly noncompact) finite convex polyhedron. A natural
question arises; how do we determine the Thurston norm on
$H^1(M;\R)$?

To address this question McMullen \cite{Mc02} used a homological
approach. It is well-known that for a knot $K$ in the 3-sphere
\[ 2\,\mbox{genus}(K) \geq \deg\left(\Delta_K(t)\right),   \]
where $\Delta_K(t) \in \Z\tpm$ denotes the Alexander polynomial of $K$. Generalizing this
McMullen \cite{Mc02} considered the multivariable Alexander polynomial $\Delta_M \in
\Z[FH_1(M;\Z)]$ (cf. Section \ref{sectiontwialex} for a definition) where
$FH_1(M;\Z):=H_1(M;\Z)/\tor_\Z(H_1(M;\Z))$ is the maximal free abelian quotient of $H_1(M;\Z)$.
Using the multivariable Alexander polynomial he defined another seminorm (called \emph{the
Alexander norm of $M$}) $||-||_A$ on $H^1(M;\R)$  as follows. If $\Delta_{M}=0$ then we set
$||\phi||_{A}=0$ for all $\phi\in H^1(M;\R)$. Otherwise for $\Delta_{M}=\sum a_if_i$ with
$a_i\in \Z$ and $f_i \in FH_1(M;\Z)$ and given $\phi \in H^1(M;\R)$ we define
\[ ||\phi||_{A} :=\sup \phi(f_i-f_j).\]
with the supremum over $(f_i, f_j)$ such that $a_ia_j\ne 0$. Note
that $\phi\in H^1(M;\R)$ naturally induces a homomorphism
$H_1(M;\R) \to \R$.

The Alexander norm ball is again a (possibly noncompact) finite convex polyhedron. McMullen showed
that the Alexander norm gives a lower bound on the Thurston norm. More precisely he proved the
following theorem.

\begin{theorem}\label{thmalexnorm}\cite[Theorem 1.1]{Mc02}
Let $M$ be a  3-manifold whose boundary is empty or consists of
tori. Then the Alexander and Thurston norms on $H^1(M;\Z)$ satisfy
\[ ||\phi||_T \geq ||\phi||_A -
\left\{ \ba{ll} 1+b_3(M), &\mbox{ if }b_1(M)=1 \mbox{ and } H^1(M;\Z) \mbox{ is generated by } \phi, \\
0, &\mbox{ if }b_1(M)>1. \ea \right.
\]
Equality holds if $\phi : \pi_1(M) \to \Z$ is represented by a fibration $M \to S^1$ such that
$M\ne S^1\times D^2$ and $M\ne S^1\times S^2$.
\end{theorem}


In \cite{Mc02}, using the Alexander norm, McMullen completely determined the Thurston norm of
many link complements. The computation was based on the following observation for the case
$b_1(M) >1$.
\\


\noindent \emph{Observation:} The Thurston norm ball lies inside
the Alexander norm ball. If the Alexander norm ball and the
Thurston norm ball agree on all extreme vertices of the Alexander
norm ball, then they agree everywhere by convexity.
\\

Note that Seiberg-Witten theory \cite{KM97} and Heegard-Floer
homology \cite{OS04} can be used to completely determine the
Thurston norm  (cf. \cite{Kr98, Kr99, Vi99,Vi03}),
but computations are not combinatorial and sometimes difficult to
apply in practice. In this paper we will take a homology theoretic
approach and find lower bounds on the Thurston norm which are easily
computed in a combinatorial way.


McMullen's homological approach has been generalized by many authors. In
\cite{Co04,Ha05,Tu02b,FK05} much stronger lower bounds for $||\phi||_T$ for \emph{specific}
$\phi \in H^1(M;\R)$ were found. In particular when $b_1(M)=1$ these methods allow us to
determine the Thurston norm ball in many cases. For the case $b_1(M) >1$ Turaev introduced
\emph{the torsion norm} generalizing McMullen's Alexander norm using \emph{abelian}
representations \cite[Chapter 4]{Tu02a}.
In this paper, given any finite dimensional representation over a
field, we define the \emph{twisted Alexander norm} and prove that it
gives a lower bound on the Thurston norm. This generalizes the work
of McMullen \cite{Mc02} and Turaev \cite{Tu02a}. Note that in a
separate paper the first author and Shelly Harvey \cite{FH06} will
show that the invariants in \cite{Ha05} are a norm as well.

In the following let $\F$ be a commutative field and
$\a:\pi_1(M)\to \glfk$ a representation. Then we define \emph{the
twisted multivariable Alexander polynomial} $\Delta_M^\a\in
\F[FH_1(M;\Z)]$ associated to $\a$ and the natural surjection
$\pi_1(M) \to FH_1(M;\Z)$ (see Section \ref{sectiontwialex}).
Similarly to the way the multivariable Alexander polynomial gives
rise to the Alexander norm we use the twisted multivariable
Alexander polynomial to define the twisted Alexander norm
$||-||_A^\a$ on $H^1(M;\R)$ associated to $\a$ (see Section
\ref{sec:twisted alexander norm}).

Let $\phi \in H^1(M;\Z)$. This defines a homomorphism
$\phi:\pi_1(M)\to \Z\cong \langle t^{\pm 1}\rangle$. We now define
$\Delta_{\phi}^{\a,i}(t) \in \ft$ to be the order of the $i$--th
twisted homology module $H_i^\a(M;\f^k \otimes_\F \ft)$ associated
to $\a$ and $\phi$. (See Section \ref{sectiontwialex}. We also refer
to \cite{KL99, FK05}.) We write $\Delta^\a_\phi(t)$ for
$\Delta_{\phi}^{\a,1}(t)$. The notion of twisted Alexander
polynomial originated from a preprint of Lin \cite{Lin01} from 1990
and was developed by Wada \cite{Wa94}. The homological definition of
twisted Alexander polynomials, which we use in this paper, was first
introduced by Kirk and Livingston \cite{KL99}. We also refer to
\cite{Kit96, FK05} for more about twisted Alexander polynomials.

In \cite[Theorem~1.1]{FK05} the authors show that twisted one-variable Alexander polynomials give
lower bounds on $||\phi||_T$ for \emph{specific} $\phi\in H^1(M;\Z)$. The following theorem allows
us to translate bounds on $||\phi||_T$ for specific $\phi\in H^1(M;\Z)$ from \cite{FK05} to bounds
on $||-||_T$ given by twisted Alexander norms. Note that $\phi$ induces a homomorphism $\phi :
\F[FH_1(M;\Z)] \to \ft$.


\newtheorem*{thm1a}{Theorem~\ref{mainthmalex}}
\begin{thm1a}
\it{ Let $M$ be a 3-manifold with $b_1(M)>1$ whose boundary is
empty or consists of tori. Let $\a : \pi_1(M) \to \glfk$ be a
representation. Let $\phi\in H^1(M;\Z)$. Then
\[ \Delta^\a_\phi(t)=\phi(\Delta^\a_M){\twialexzero}\twialextwo.\]
Furthermore if $\phi(\Delta^\a_M)\ne 0$, then  $\twialexzero\ne 0$
and $\twialextwo\ne 0$ and hence $\twialex\ne 0$.}
\end{thm1a}

The proof is based on the functoriality of Reidemeister torsion (see Section \ref{sec:multi})
and builds on ideas of Turaev. The following two theorems are our main results.

\newtheorem*{thm1}{Theorem~\ref{mainthm}}
\begin{thm1}[{\bf Main Theorem 1}]
\it{ Let $M$ be a 3-manifold with $b_1(M)>1$ whose boundary is
empty or consists of tori. Let $\a:\pi_1(M)\to \glfk$ be a
representation. Then for the corresponding twisted Alexander norm
$||-||^\a_A$, we have
\[ ||\phi||_T \ge \frac1k ||\phi||^\a_A   \]
for all $\phi\in H^1(M;\R)$. }
\end{thm1}

Let $M$ be a 3-manifold and $\phi\in H^1(M;\Z)$. We say
\emph{$(M,\phi)$ fibers over $S^{1}$} if the homotopy class of maps
$M\to S^1$ induced by $\phi:\pi_1(M)\to H_1(M;\Z)\to \Z$ contains a
representative that is a fiber bundle over $S^{1}$. Thurston
\cite{Th86} showed that if $(M,\phi)$ fibers over $S^1$, then $\phi$
lies in the cone  on a top-dimensional open face of the Thurston
norm ball. We denote this cone by $C(\phi)$.

\newtheorem*{thm2}{Theorem~\ref{mainthmfib}}
\begin{thm2}[{\bf Main Theorem 2}]
\it{ Let $M$ be a   3-manifold with $b_1(M)>1$ whose boundary is empty or consists of tori such
that $M\ne S^1\times D^2$ and $M\ne S^1\times S^2$. Let $\a:\pi_1(M)\to \glfk$ be a
representation. If $\phi\in H^1(M;\Z)$ is such that $(M,\phi)$ fibers over $S^1$, then
\[ ||\psi||_T = \frac1k ||\psi||^\a_A   \]
 for all $\psi\in C(\phi)$.
}
\end{thm2}

By Theorem \ref{mainthm} twisted Alexander norms give lower bounds on the Thurston norm. With
the same reason as for the Alexander norm ball, twisted Alexander norm balls are (possibly
noncompact) finite convex polyhedra. Therefore we can use McMullen's observation in the above to
determine the Thurston norm using twisted Alexander norms.

\begin{figure}[h] \begin{center}
\includegraphics[scale=0.25]{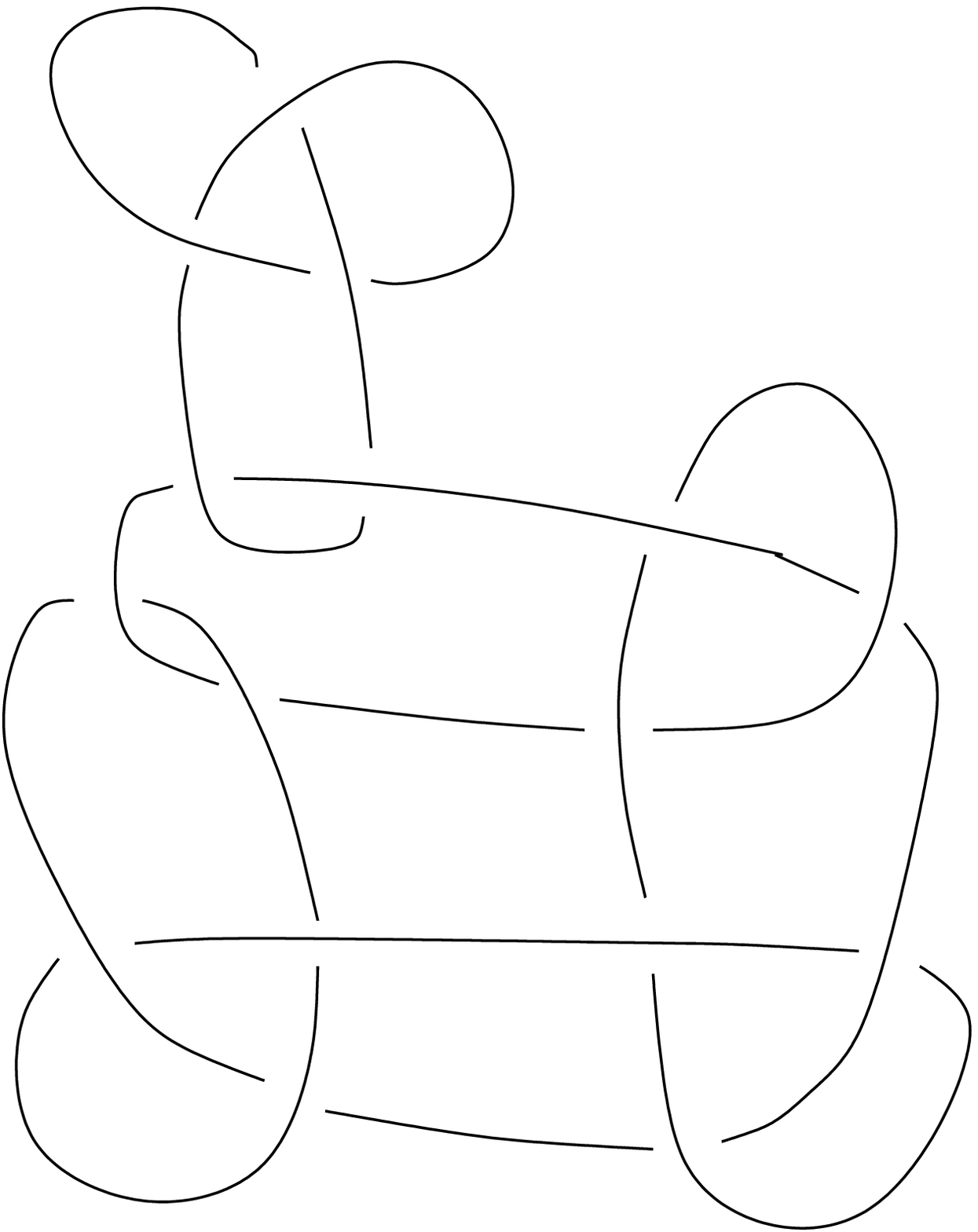}
\caption{Link $L$} \label{link11n73intro}
 \end{center}
 \end{figure}

\begin{figure}[h] \begin{center}
\begin{tabular}{rl}
\includegraphics[scale=0.3]{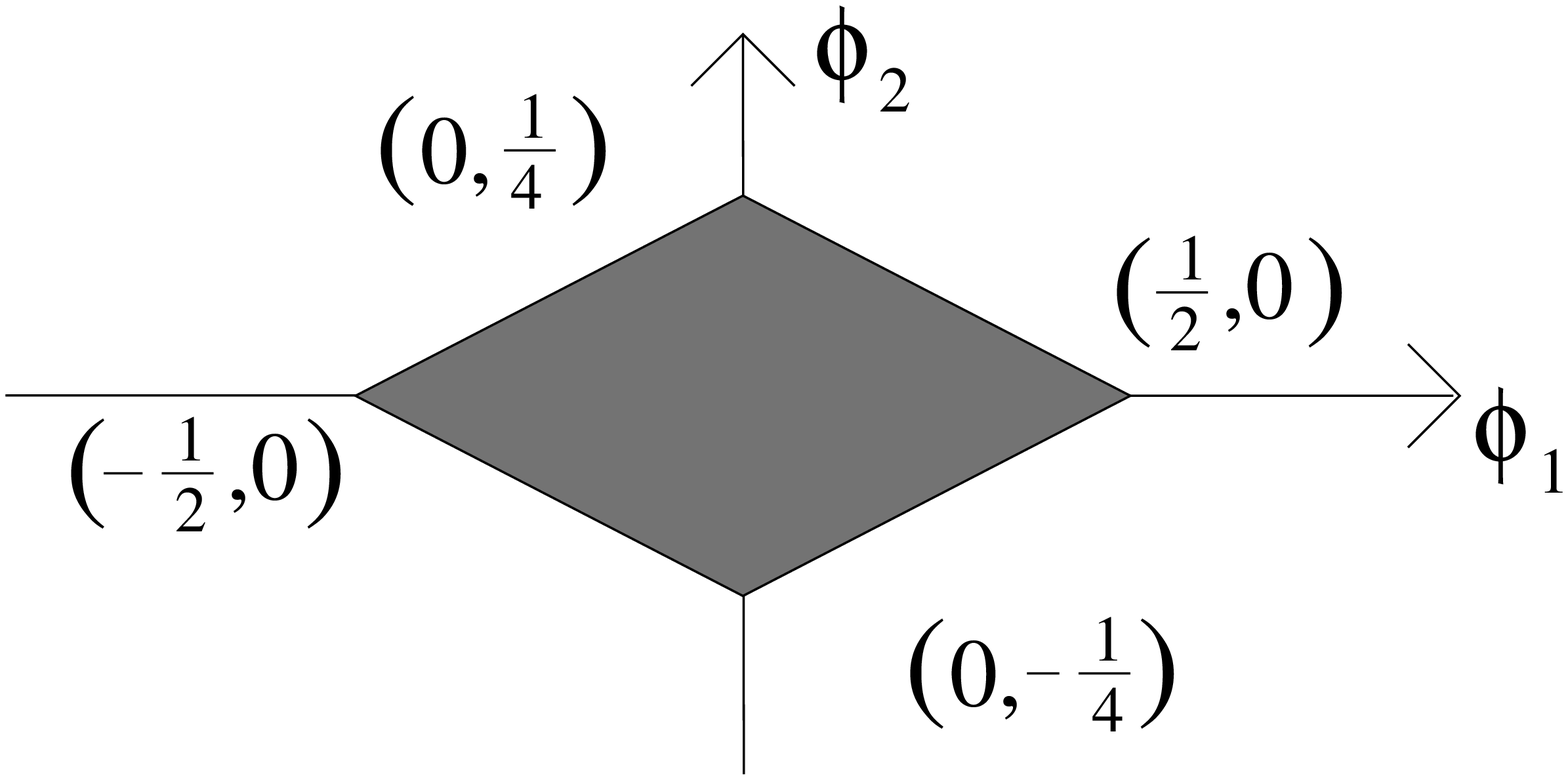}&\includegraphics[scale=0.3]{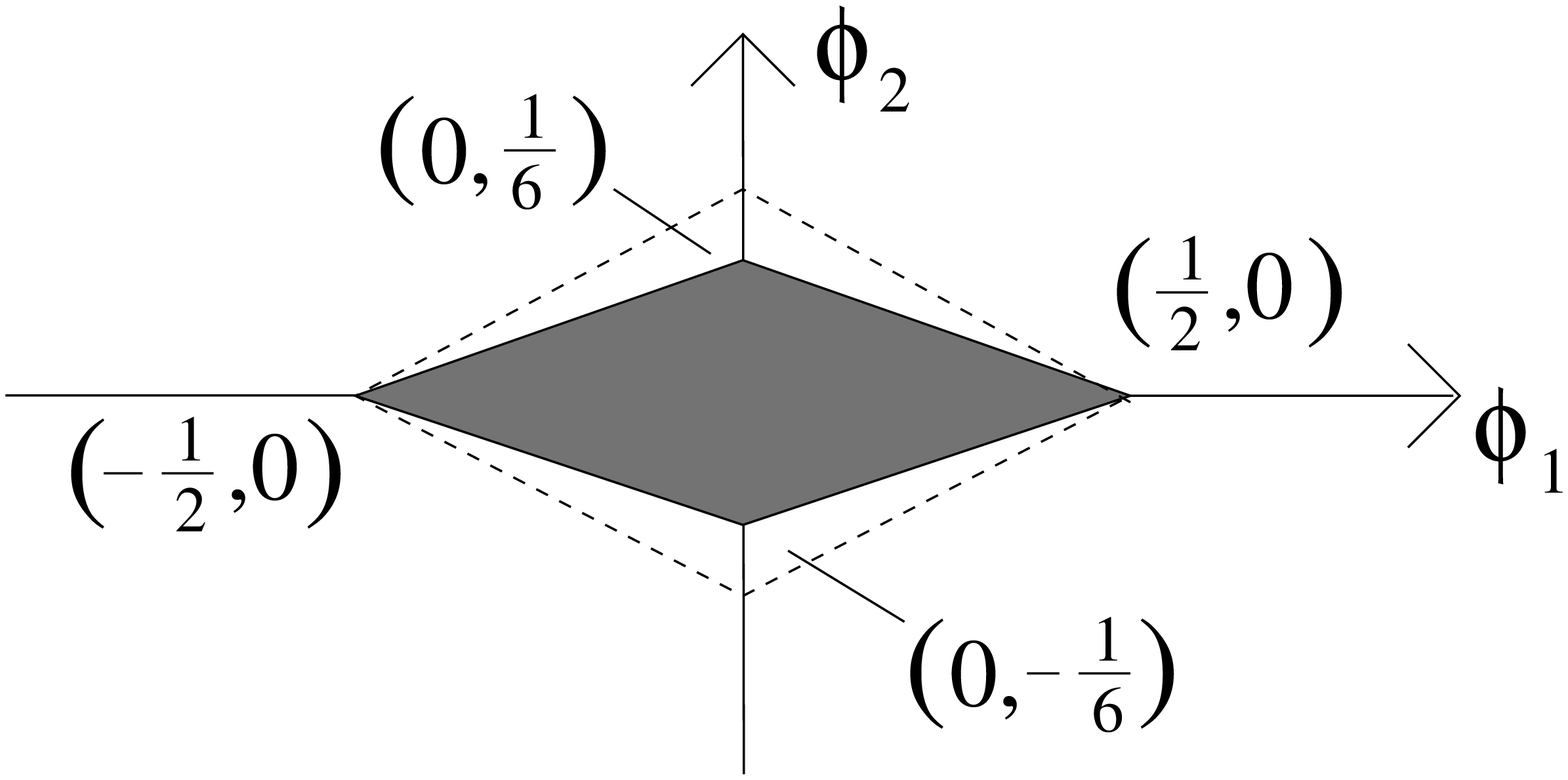}
\end{tabular}
 \caption{The untwisted and the twisted Alexander norm ball of $L$.} \label{normballintro}
\end{center} \end{figure}

In Section \ref{sec:example} we give examples which show how
powerful twisted Alexander norms are. For example we determine the
Thurston norm of the complement of the link $L$ in Figure
\ref{link11n73intro}, which can not be determined by the (usual)
Alexander norm. The components of $L$ are $K_1$, the trefoil, and
$K_2=11_{440}$ (here we use \emph{knotscape} notation).
Let $X(L)$ denote the complement of an open tubular neighborhood of $L$ in the 3-sphere. Then
\[ \Delta_{X(L)}(x_1,x_2)=(x_1^2-x_1+1)(x_2^4-2x_2^3+3x_2^2-2x_2+1)\in \Q
\xypm.\] The resulting
Alexander norm ball is given in Figure \ref{normballintro} on the left. On the other hand using the
program \emph{KnotTwister} \cite{F05} we found a representation $\a:\pi_1(X(L))\to \gl(\F_{13},2)$
such that
\[ \Delta_{X(L)}^\a(x_1,x_2)=\Delta_1(x_1)\Delta_2(x_2)\]
where $\deg(\Delta_1(x_1))=4$ and $\deg(\Delta_2(x_2))=12$. (Here
$\F_n$ denotes the field of $n$ elements.) Hence the twisted
Alexander norm ball for $\frac{1}{2}||-||_A^\a$ is the shaded
region given in Figure \ref{normballintro} on the right.
By Theorem \ref{mainthm} we have $||\phi||_T \ge \frac{1}{2}
||\phi||_{A}^\a$. It is clear from Figure \ref{normballintro} that
$\frac12||-||_A^\a$ gives a strictly sharper bound on the Thurston
norm than $||-||_A$ does. In Section \ref{sec:hophlike} we will
see that the norms $||-||_T$ and $\frac12||-||_A^\a$ agree on the
vertices of the norm ball of $\frac12||-||_A^\a$. Therefore by
McMullen's observation the norms agree everywhere.
Hence the shaded region in Figure \ref{normballintro} on the right
is in fact the Thurston norm ball of the link $L$. We point out
that it follows immediately from Theorem \ref{mainthmfib} that
$(X(L),\phi)$ does not fiber over $S^1$ for any $\phi \in
H^1(M;\Z)$. See Section \ref{sec:example} for more details.



Our approach works very well in many cases, but sometimes it is
difficult to find an appropriate representation. Therefore it is
sometimes convenient to find lower bounds on the Thurston norm of
a finite cover $\ti{M}$ of $M$. By a result of Gabai
\cite[p.~484]{Ga83} (cf. also Theorem \ref{lemmathurstong}) the
Thurston norm on $\ti{M}$ determines the Thurston norm on $M$. In
many cases it is easier to find representations of $\ti{M}$. This
approach allows us to determine the Thurston norm ball of
Dunfield's link \cite{Du01} (see Section \ref{sec:dunfield}).
\\

{\bf Outline of the paper:} In Section \ref{sec:twisted
invariants} we define twisted Alexander modules and twisted
Alexander polynomials. In Section \ref{sec:bound} we define
twisted Alexander norms and prove the main theorems. We quickly
discuss how to compute twisted Alexander polynomials in Section
\ref{sec:computation} and give examples in Section
\ref{sec:example}. In Section \ref{sec:multi} we give a proof of
Theorem \ref{mainthmalex} which shows the precise relationship
between the twisted multivariable Alexander polynomials and the
twisted one-variable Alexander polynomials.
\\

{\bf Notations and conventions:} For a link $L$ in $S^3$, $X(L)$
denotes the exterior of $L$ in $S^3$. (That is, $X(L) = S^3\sm \nu
L$ where $\nu L$ is an open tubular neighborhood of $L$ in $S^3$.)
An arbitrary (commutative) field is denoted by $\F$. $\F_n$
denotes the finite field of $n$ elements. We identify the group
ring $\F[\Z]$ with $\ft$. We denote the permutation group of order
$k$ by $S_k$. For a 3-manifold $M$ we use the canonical
isomorphisms to identify $H^1(M;\Z) = \Hom(H_1(M;\Z), \Z) =
\Hom(\pi_1(M),\Z)$. Hence sometimes $\phi\in H^1(M;\Z)$ is
regarded as a homomorphism $\phi : \pi_1(M) \to \Z$ (or $\phi :
H_1(M;\Z) \to \Z$) depending on the context.
\\

{\bf Acknowledgments:} The authors would like to thank Stefano
Vidussi and Jae Choon Cha for helpful conversations and
suggestions.

\section{Twisted Alexander polynomials}
\label{sec:twisted invariants}

In this section we give the definition of twisted Alexander
polynomials.
\subsection{Torsion invariants}
Let $R$ be a commutative Noetherian unique factorization domain
(henceforth UFD). An example of $R$ to keep in mind is
$\F[t_1^\pm, t_2^\pm, \ldots, t_n^\pm]$, a (multivariable) Laurent
polynomial ring over a field $\F$. For a finitely generated
$R$-module $A$, we can find a presentation
$$
R^r \xrightarrow{P} R^s \to A \to 0
$$
since $R$ is Noetherian. Let $i\ge 0$ and suppose $s-i\le r$. We
define $E_i(A)$, \emph{the $i$-th elementary ideal} of $A$, to be
the ideal in $R$ generated by all $(s-i)\times (s-i)$ minors of
$P$ if $s-i>0$ and to be $R$ if $s-i\le 0$. If $s-i > r$, we
define $E_i(A)= 0$. It is known that $E_i(A)$ does not depend on
the choice of a presentation of $A$ (cf. \cite{CF77}).

Since $R$ is a UFD there exists a unique smallest principal ideal
of $R$ that contains $E_0(A)$. A generator of this principal ideal
is defined to be the \emph{order of $A$} and denoted by $\ord
(A)\in R$. The order is well-defined up to multiplication by a
unit in $R$. Note that  $A$ is not $R$-torsion if and only if
$\ord (A) =0$. For more details, we refer to \cite{Hi02}.

\subsection{Twisted Alexander invariants}\label{sectiontwialex}
Let $M$ be a 3-manifold and $\psi:\pi_1(M)\to F$ a homomorphism to
a free abelian group $F$. We do not demand that $\psi$ is
surjective. Note that $\L:=\F[F]$ is a commutative Noetherian UFD.
Let $\a:\pi_1(M)\to \gl(\F,k)$ be a representation.

Using $\a$ and $\psi$, we define a left $\Z[\pi_1(M)]$-module
structure on $\F^k\otimes_\F \L=:\L^k$ as follows:
\[  g\cdot (v\otimes p):= (\a(g)\cdot v)\otimes (\psi(g)p)\]
where $g\in \pi_1(M)$ and $v\otimes p \in \F^k\otimes_\F \L=\L^k$.
Together with the natural structure of $\L^k$ as a $\L$--module we
can view $\L^k$ as a $\Z[\pi_1(M)]$--$\L$ bi--module.
Recall there exists a canonical left $\pi_1(M)$--action on the universal cover $\ti{M}$. We
consider the chain complex $C_*(\ti{M})$ as a right $\Z[\pi_1(M)]$-module by defining $\s \cdot
g:=g^{-1}\s$ for a singular chain $\s$. For $i\ge 0$, we define \emph{the $i$-th twisted
Alexander module of $(M,\psi,\a)$} to be
$$
H_i^\a(M;\L^k) := H_i(C_*(\tilde{M})\otimes_{\Z[\pi_1(M)]}\L^k).
$$
Since $\L^k$ is a right $\L$-module twisted Alexander modules can
be regarded as right $\L$-modules. Since $M$ is compact and $\L$
is Noetherian these modules are finitely generated over $\L$.

\begin{defn} \label{def:polynomial}
The \emph{$i$-th (twisted) Alexander polynomial of $(M,\psi,\a)$}
is defined to be $\ord (H_i^\a(M;\L^k))\in \L$ and denoted by
$\Delta^{\a,i}_{M,\psi}$. When $i=1$, we drop the superscript $i$
and abbreviate $\Delta^{\a,i}_{M,\psi}$ by $\Delta^\a_{M,\psi}$,
and we call it the \emph{(twisted) Alexander polynomial of
$(M,\psi,\a)$}.
\end{defn}

Twisted Alexander polynomials are well-defined up to
multiplication by a unit in $\L$. We drop the notation $\psi$ when
$\psi$ is the natural surjection to $FH_1(M;\Z)$. We also drop
$\a$ when $\a$ is the trivial representation to $\gl(\Q,1)$ and
drop $M$ in the case that $M$ is clear from the context. If $\psi$
is a homomorphism to $\Z$ then we identify $\F[\Z]$ with $\ft$ and
we write $\Delta_{M,\psi}^{\a,i}(t)\in \ft$. The above homological
definition of twisted Alexander polynomials was first introduced
by Kirk and Livingston \cite{KL99}.

\section{Twisted Alexander norms as lower bounds on the Thurston norm}
\label{sec:bound}

In this section we define twisted Alexander norms, which
generalize the Alexander norm of McMullen \cite{Mc02} and the
torsion norm of Turaev \cite{Tu02a}. We show that twisted
Alexander norms give lower bounds on the Thurston norm and that
they give fibering obstructions of 3-manifolds.
\subsection{Twisted Alexander norm} \label{sec:twisted alexander
norm} Following an idea of McMullen's \cite{Mc02} we now use the twisted multivariable Alexander
polynomial corresponding to $\psi:\pi_1(M)\to FH_1(M;\Z)$ to define a norm on $H^1(M;\R)$. Let
$\a:\pi_1(M)\to \gl (\F,k)$ be a representation. If $\Delta_{M}^{\a}=0$ then we set
$||\phi||_{A}^{\a}=0$ for all $\phi\in H^1(M;\R)$. Otherwise we write $\Delta_{M}^{\a}=\sum a_if_i$
for $a_i\in \F$ and $f_i \in FH_1(M;\Z)$. Given $\phi \in H^1(M;\R)$ we then define
\[ ||\phi||_{A}^{\a} :=\sup \phi(f_i-f_j),\]
with the supremum over $(f_i, f_j)$ such that $a_ia_j\ne 0$.
Clearly this defines a seminorm on $H^1(M;\R)$ which we call the
\emph{twisted Alexander norm of $(M,\a)$}. This is a
generalization of the Alexander norm introduced by McMullen
\cite{Mc02}. Indeed, the Alexander norm is the same as the twisted
Alexander norm corresponding to the trivial representation $\a :
\pi_1(M) \to \gl (\Q,1)$. In this case we just write $||-||_{A}$.
Twisted Alexander norms also generalize the torsion norm of Turaev
\cite{Tu02a}.
\subsection{Lower bounds on the Thurston norm}

Recall that McMullen showed that in the case $b_1(M)>1$ the
Alexander norm $||-||_A$ is a lower bound on the Thurston norm
(see Theorem \ref{thmalexnorm}). We extend this result to twisted
Alexander norms.

%

\begin{theorem}[{\bf Main Theorem 1}]\label{mainthm}
Let $M$ be a 3-manifold with $b_1(M)>1$ whose boundary is empty or
consists of tori. Let $\a:\pi_1(M)\to \glfk$ be a representation.
Then for the corresponding twisted Alexander norm $||-||^\a_A$ we
have
\[ ||\phi||_T \ge \frac1k ||\phi||^\a_A   \]
for all $\phi\in H^1(M;\R)$.
\end{theorem}
\noindent This theorem generalizes McMullen's theorem (Theorem
\ref{thmalexnorm}). Turaev \cite{Tu02a} proved this theorem in the
special case of abelian representations.


\begin{theorem}[{\bf Main Theorem 2}]\label{mainthmfib}
Let $M$ be a 3-manifold with $b_1(M)>1$ whose boundary is empty or
consists of tori such that $M\ne S^1\times D^2$ and $M\ne
S^1\times S^2$ . Let $\a:\pi_1(M)\to \glfk$ be a  representation.
If $\phi\in H^1(M;\Z)$ is such that $(M,\phi)$ fibers over $S^1$,
then
\[ ||\psi||_T = \frac1k ||\psi||^\a_A   \]
 for all $\psi\in C(\phi)$.
\end{theorem}



\noindent The idea of the proofs of the main theorems is to combine
the lower bounds for one-variable Alexander polynomials from
\cite{FK05} with Theorem \ref{mainthmalex}. In \cite{FK05} we proved
the following theorem.


\begin{theorem} \cite[Theorem~1.1~and~Theorem~1.2]{FK05} \label{mainthmfk05}
Let $M$ be a 3-manifold whose boundary is empty or consists of
tori. Let $\phi \in H^1(M)$ be nontrivial and $\a:\pi_1(M)\to
\glfk$ a representation such that $\twialex \ne 0$. Then
$\twialexi\ne 0$ for $i=0,2$ and
\[ \tnphi \geq
\frac{1}{k}\big(\degtwialex- \degtwialexzero -\degtwialextwo \big).\] Furthermore, if $(M,\phi)$
fibers over $S^1$ and if $M\ne S^1\times D^2$ and $M\ne S^1\times S^2$, then equality holds.
\end{theorem}

We also need the following theorem to prove the main theorems.
This theorem clarifies the precise relationship between the
twisted multivariable Alexander polynomial and the twisted
one-variable Alexander polynomials of a 3-manifold.

\begin{theorem}\label{mainthmalex}
Let $M$ be a 3-manifold with $b_1(M)>1$ whose boundary is empty or
consists of tori. Let $\a : \pi_1(M) \to \glfk$ be a
representation. Let $\phi\in H^1(M;\Z)$ be nontrivial. Then
\[ \twialex=\phi(\Delta^\a_M){\twialexzero}\twialextwo. \]
Furthermore if $\phi(\Delta^\a_M)\ne 0$, then  $\twialexzero\ne 0$ and $\twialextwo\ne 0$ and hence
$\twialex\ne 0$.
\end{theorem}

\noindent The idea of the proof of Theorem \ref{mainthmalex} is to
go from the twisted multivariable Alexander polynomials to
Reidemeister torsion which is functorial, and then to go back to
the twisted one-variable Alexander polynomials. The proof of
Theorem \ref{mainthmalex} is postponed to Section
\ref{sec:funtoriality}. Now we give a proof of Theorem
\ref{mainthm}.

\begin{proof}[Proof of Theorem \ref{mainthm}]
If $\Delta_M^\a = 0$, then $||\phi||^\a_A = 0$ for all $\phi\in
H^1(M;\R)$, hence the theorem holds. We now consider the case
$\Delta_M^\a \ne 0$.

First suppose that $\phi\in H^1(M;\Z)$ is nontrivial and lies
inside the cone on an open top-dimensional face of the twisted
Alexander norm ball. Write $\Delta_M^\a = \sum a_i f_i$ where $a_i
\in \F \setminus \{0\}$ and $f_i \in FH_1(M;\Z)$. We have
\[
\phi(\Delta_M^\a) = \sum a_i t^{\phi(f_i)}
\]
in $\ft$. Since $\phi$ is inside the cone on an open
top-dimensional face of the twisted Alexander norm ball, the
highest and lowest values of $\phi(f_i)$ occur only once in the
above equation. Therefore $\phi(\Delta_M^\a) \ne 0$ and
\[
\deg \left(\phi(\Delta_M^\a)\right) =||\phi||^\a_A.
\]
By Theorem \ref{mainthmalex} we have $\twialex \ne 0$, $\twialexzero \ne 0$, $\twialextwo\ne 0$ and
\be \label{equn1} \deg \left(\Delta_\phi^\a(t) \right) =
||\phi||^\a_A+\deg\left(\twialexzero\right)+ \deg\left(\twialextwo\right). \ee
 Since $\twialex\ne 0$ we get by Theorem
\ref{mainthmfk05} that
 \be \label{equn2} \tnphi \geq \frac{1}{k}\left( \degtwialex-
\deg\left(\twialexzero\right)-\deg\left(\twialextwo\right)\right). \ee
 Combining the inequalities
(\ref{equn1}) and (\ref{equn2}) we clearly get $\tnphi \geq \frac{1}{k} ||\phi||^\a_A$.
This proves Theorem \ref{mainthm} for all  $\phi\in H^1(M;\Z)$
inside the cone on an open top-dimensional face of the twisted
Alexander norm ball. By homogeneity and continuity we get that in
fact $||\phi||_T \ge \frac{1}{k} ||\phi||^\a_A$ for all $\phi\in
H^1(M;\R)$.
\end{proof}

For the proof of Theorem \ref{mainthmfib} we need the following theorem proved by Thurston
\cite{Th86} and which can also be found in \cite[Theorem 9, p.~259]{Oe86}.
\begin{theorem}[Thurston]
\label{fiberface} Let $M$ be a   3-manifold.  If $\phi\in H^1(M;\Z)$ is such that $(M,\phi)$ fibers
over $S^1$, then $\phi$ lies in the cone on a top-dimensional open face of the Thurston norm ball.
Furthermore, if we denote this cone by $C(\psi)$, then $(M,\psi)$ fibers over $S^1$ for all $\psi
\in C(\psi)\cap H^1(M;\Z)$.
\end{theorem}

\begin{proof}[Proof of Theorem \ref{mainthmfib}]
Suppose $\phi\in H^1(M;\Z)$. If $\phi$ is nontrivial and $(M,\phi)$ fibers over $S^1$ then the
inequality in Theorem \ref{mainthmfk05} and hence by the proof of Theorem \ref{mainthm} the
inequality in Theorem \ref{mainthm} become equalities. Furthermore, by  Theorem \ref{fiberface},
$\phi$ lies in the cone on a top-dimensional open face $C(\psi)$ of the Thurston norm ball, and
$(M,\psi)$ fibers over $S^1$ for any $\psi\in C(\psi)\cap H^1(M;\Z)$. In particular  we have
\[
||\psi||_T = \frac1k
||\psi||^\a_A
\] for every $\psi\in C\cap H^1(M;\Z)$ which is
nontrivial. By homogeneity and continuity it follows that
\[
||\psi||_T=\frac{1}{k}||\psi||_A^{\a}
\]
for all $\psi\in C(\psi)$.
\end{proof}

\section{Computation of twisted Alexander norms}
\label{sec:computation}

Let $M$ be a  3-manifold and $\psi:\pi_1(M)\to F$ a homomorphism
to a free abelian group $F$ such that $\psi:H_1(M;\Q)\to
F\otimes_\Z \Q$ is surjective. (In this case we say $\psi$ is
\emph{rationally surjective}.) Given a representation
$\a:\pi_1(M)\to \glfk$ we quickly outline how to compute
$\Delta_{M,\psi}^{\a}$ and hence the twisted Alexander norm.

Denote the universal cover of $M$ by $\ti{M}$. If $p$ is a point in $M$, then denote the preimage
of $p$ under the map $\ti{M}\to M$ by $\ti{p}$. Then a presentation matrix for
$$
H_i^\a(M, p;\F^k[F]):=H_i(C_*(\tilde{M},\ti{p})\otimes_{\Z[\pi_1(M)]}\F^k[F]).
$$
can be found using Fox calculus from a presentation of the group $\pi_1(M)$.
We also refer to the literature \cite{Fo53,Fo54,CF77}), but we point out that we view $C_*(\ti{M})$
as a \emph{right} $\Z[\pi_1(M)]$-module, whereas the literature normally views $C_*(\ti{M})$ as a
\emph{left} $\Z[\pi_1(M)]$-module  (cf. also \cite[Section~6]{Ha05}).

By using the long exact sequence of the twisted homology modules
of the pair of spaces $(M,p)$, one can obtain the following short
exact sequence of $\F[F]$-modules:
 $$ 0\to H_1^\a(M;\F^k[F]) \to
H_1^\a(M,p;\F^k[F]) \to A \to 0 $$
 where $A=\ker\{H_0^\a(p;\F^k[F]) \to H_0^\a(M;\F^k[F])\}$. Note
that $H_0^\a(p;\F^k[F])\cong \F^k[F]$ whereas $H_0^\a(M;\F^k[F])$ is
a finite-dimensional $\F$-vector space by the following well--known
lemma.
\begin{lemma} \label{lemmah0m}
Let $X$ be a 3-manifold, $\psi:\pi_1(X) \to F$ a rationally
surjective map with $F$ a free abelian group, and $\a:\pi_1(X)\to
\glfk$ a representation. Then
\[ H_i^\a(X;\fk[F])=H_i(\ker(\psi);\fk)^n, i=0,1\]
where $n = |F/\im(\psi)|$.
\end{lemma}

It follows that  $A$ is an $\F[F]$-module of rank $k$. (For the
notion of rank over $\F[F]$ we refer to the first paragraph in
Section \ref{sec:calc}.) If $H_0^\a(M;\fk[F])$ is $\F[F]$-torsion,
then by \cite[Theorem 3.4]{Hi02}
$$
\Delta^\a_{M,\psi} = \ord(E_0(H_1^\a(M;\fk[F])))= \ord
(E_k(H_1^\a(M,p;\fk[F]))),
$$
which can be computed using the presentation matrix for
$H_1^{\a}(M,p;\F^k[F])$. If $H_1^\a(M;\F^k[F])$ is not
$\F[F]$-torsion, $E_k(H_1^\a(M,p;\fk[F])) = 0$ and
$\Delta^\a_{M,\psi} = \ord (E_k(H_1^\a(M,p;\fk[F]))) = 0$.

In the case that $\partial M\ne 0$  we can compute
$\Delta_{M,\psi}^{\a}$ from Wada's invariant, which tends to be
easier to compute.  We refer to \cite{Wa94, KL99} for more
details.
\section{Examples for twisted Alexander norms}
\label{sec:example}

In this section, using twisted Alexander norms, we completely
determine the Thurston norm of two examples: certain Hopf-like
links and Dunfield's link \cite{Du01}.

\subsection{Hopf-like links} \label{sec:hophlike}
\begin{figure}[h]
\begin{center}
\begin{tabular}{rl}
\includegraphics[scale=0.25]{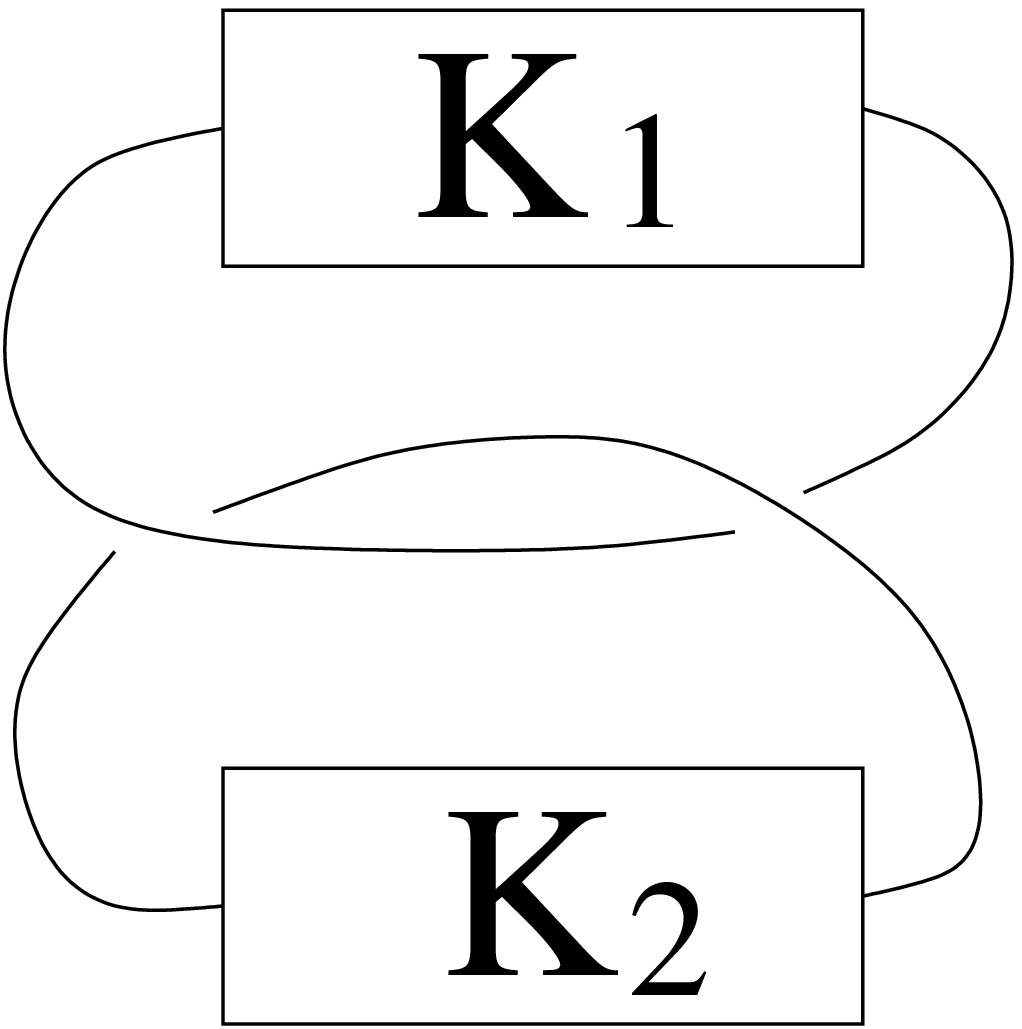}& \hspace{1cm}
\includegraphics[scale=0.25]{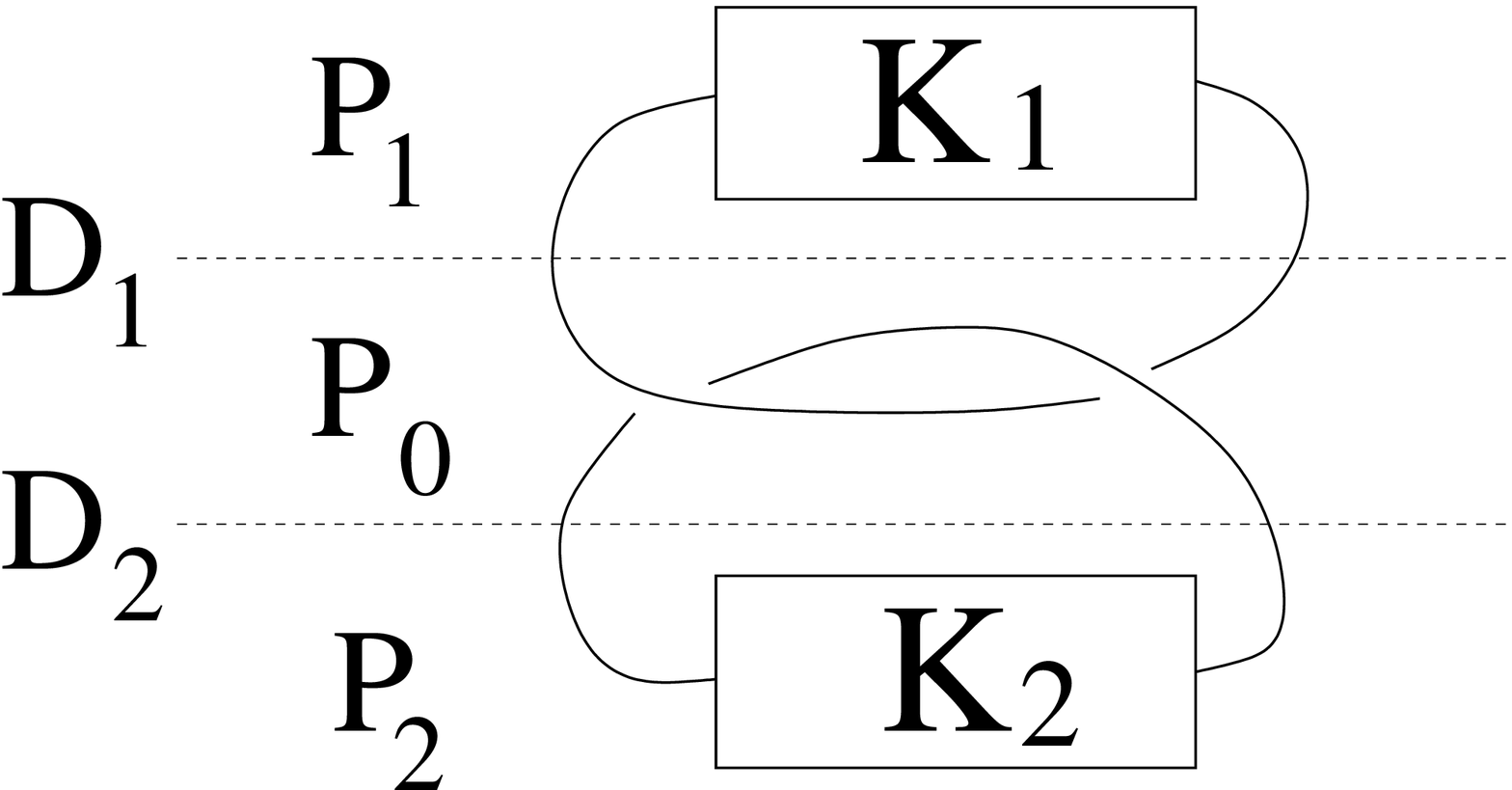}
\end{tabular}
 \caption{The link $L$ and the link complement
cut along annuli $D_1$ and $D_2$} \label{linkcut}\label{hopflink}\label{linkk1k2}
\end{center}
\end{figure}

In this section, for a link $L$ (possibly with one component), we
write $\Delta_L^\a$ for $\Delta_{X(L)}^\a$. Consider a link $L$ as
in Figure \ref{linkk1k2}. We will call these links
\emph{Hopf-like}. Denote the meridian of $K_1$ by $\mu_1$ and the
meridian of $K_2$ by $\mu_2$. Denote the corresponding elements in
$H_1(X(L);\Z)$ by $x_1$ and $x_2$. We then identify
$\Z[H_1(X(L);\Z)]$ with $\Z \xypm$.

Let $D_1$ (respectively, $D_2$) be the annulus cutting through $L$ just below $K_1$
(respectively, above $K_2$). Denote the three components of $X(L)$ cut along $D_1\cup D_2$ by
$P_1, P_0, P_2$ (see Figure \ref{linkcut} below). Note that $P_i \cong X(K_i)$, $i=1,2$. In
particular any representation $\a:\pi_1(X(L))\to \glfk$ induces  representations
$\pi_1(X(K_i))\to \glfk, i=1,2,$ which we also denote by $\a$.

\begin{proposition}\label{propalexlink}
Let $\a:\pi_1(X(L))\to \glfk$ be a representation. Assume
$\Delta_{K_i}^\a(x_i)\ne 0$ for $i=1,2$. Then
\[
\Delta_L^\a(x_1,x_2)=\Delta_1(x_1)\Delta_2(x_2)\in \F \xypm
\]
where
\[
\Delta_i(x_i)=\Delta_{K_i}^\a(x_i)\frac{\det(\a(\mu_i)x_i-\id)}{\Delta_{K_i}^{\a,0}} \in \F[x_i^{\pm 1}],
\phantom{a}i=1,2.
\]
In particular
\[
\deg(\Delta_i(x_i))=\deg\left(\Delta_{K_i}^\a(x_i)\right)+k-\deg\left(\Delta_{K_i}^{\a,0}(x_i)\right),
\phantom{a} i=1,2.
\]
\end{proposition}


\begin{proof}
First note that $D_i$ is homotopy equivalent to the circle for
$i=1,2$, hence it follows from Lemma \ref{lemmah0m} that
$H_1^\a(D_i;\fk\xypm)= 0$. We now consider the Mayer-Vietoris
sequence of $X(L)=P_1\cup_{D_1}P_0\cup_{D_2}P_2$.
\[
\ba{rcccccccccccccc} 0&\hspace{-0.1cm} \to \hspace{-0.1cm} &
\bigoplus\limits_{i=0}^2 H_1^\a(P_i;\fk \xypm)&\hspace{-0.1cm}
\to\hspace{-0.1cm}
& H_1^\a(X(L);\fk \xypm)&\hspace{-0.1cm} \to\hspace{-0.1cm} &\\
\bigoplus\limits_{i=1}^2 H_0^\a(D_i;\fk \xypm)&\hspace{-0.1cm}
\to\hspace{-0.1cm} &\bigoplus\limits_{i=0}^2 H_0^\a(P_i;\fk \xypm)
&\hspace{-0.1cm} \to\hspace{-0.1cm} & H_0^\a(X(L);\fk
\xypm)&\hspace{-0.1cm} \to\hspace{-0.1cm} & 0.\ea \]
 By  \cite[Lemma 5, p.~76]{Le67} for any exact sequence of $\f
\xypm$--torsion modules the alternating product of the respective orders in $\f \xypm$ equals
one. The proposition now follows immediately from the following computations.

By  Lemma \ref{lem:delta03} we have that $\ord(H_0^\a(X(L);\fk
\xypm))=1$. We compute the orders of the twisted Alexander modules of
$P_1$ and $P_2$. Since $P_i\cong X(K_i),$ $i=1,2$, the natural
surjection $\psi : \Z[\pi_1(X(L))]\to \Z \xypm$ restricted to
$P_i$ only has values in $\Z \xipm$. Thus we get
\[ \ba{rcl}
H_j^{\a}(P_1;\fk \xypm) &\cong &H_j^\a(X(K_1);\fk \xpm)\otimes_{\f} \f \ypm \mbox{ for all } j, \mbox{ and }\\
H_j^{\a}(P_2;\fk \xypm) &\cong &H_j^\a(X(K_2);\fk
\ypm)\otimes_{\f} \f \xpm \mbox{ for all } j. \ea
\]
Therefore
\[
\ord\left(H_j^{\a}(P_i;\fk \xypm )\right) =
\Delta^{\a,j}_{K_i}(x_i).
\]
for all $j\ge 0$ and $i=1,2$.

Let us consider $P_0$. $P_0$ is homotopy equivalent to the torus and $\pi_1(P_0)$ is the free
abelian group spanned by $\mu_1$ and $\mu_2$. By  Lemma \ref{lemmah0m} we have $H_1^\a(P_0;\fk
\xypm)=0$. Therefore $\ord(H_1^\a(P_0;\fk \xypm))=1$. Furthermore the argument in the proof of
Lemma \ref{lem:delta03}  shows that $\ord(H_0^\a(P_0;\fk \xypm))=1$.

Now consider $D_1$ and $D_2$.
Using the cellular chain complex of the circle, one easily sees that
\[ \ord(H_0^\a(D_i;\fk \xypm))=\det(\a(\mu_i)x_i-\id) \]
for $i=1,2$.
\end{proof}

%
%

\begin{corollary} \label{coralexlink}
For the trivial representation $\a : \pi_1(X(L)) \to GL(\F,1)$,
\[
\Delta^\a_L(x_1,x_2)=\Delta^\a_{K_1}(x_1) \Delta^\a_{K_2}(x_2).
\]
\end{corollary}
\begin{proof}
Since $\a$ is a one-dimensional trivial representation,
\[
H_0^\a(X(K_1);\f\xpm)=\f\xpm/(x_1-1).
\]
Hence $\Delta^{\a,0}_{K_1}(x_1)=x_1-1$. Also
$\det(\a(\mu_1)x_1-\id) = x_1 - 1 = \Delta^{\a,0}_{K_1}(x_1)$.
Similarly $\det(\a(\mu_2)x_2-\id) =
\Delta^{\a,0}_{K_2}(x_2)=x_2-1$. Now use Proposition
\ref{propalexlink}.
\end{proof}

\begin{corollary} \label{corvertex}
Let $d_i:=\deg(\Delta^\a_i(x_i)),$ $i=1,2$ in Proposition
\ref{propalexlink}. Then the norm ball of $\frac1k||-||^\a_A$ has
exactly four extreme vertices namely $(\pm \frac{k}{d_1},0)$ and
$(0,\pm \frac{k}{d_2})$.
\end{corollary}
\noindent The above corollary easily follows from Proposition
\ref{propalexlink}.

\begin{figure}[h] \begin{center} \begin{tabular}{cc}
\includegraphics[scale=0.25]{link11n73.eps}&\hspace{1cm}
\includegraphics[scale=0.25]{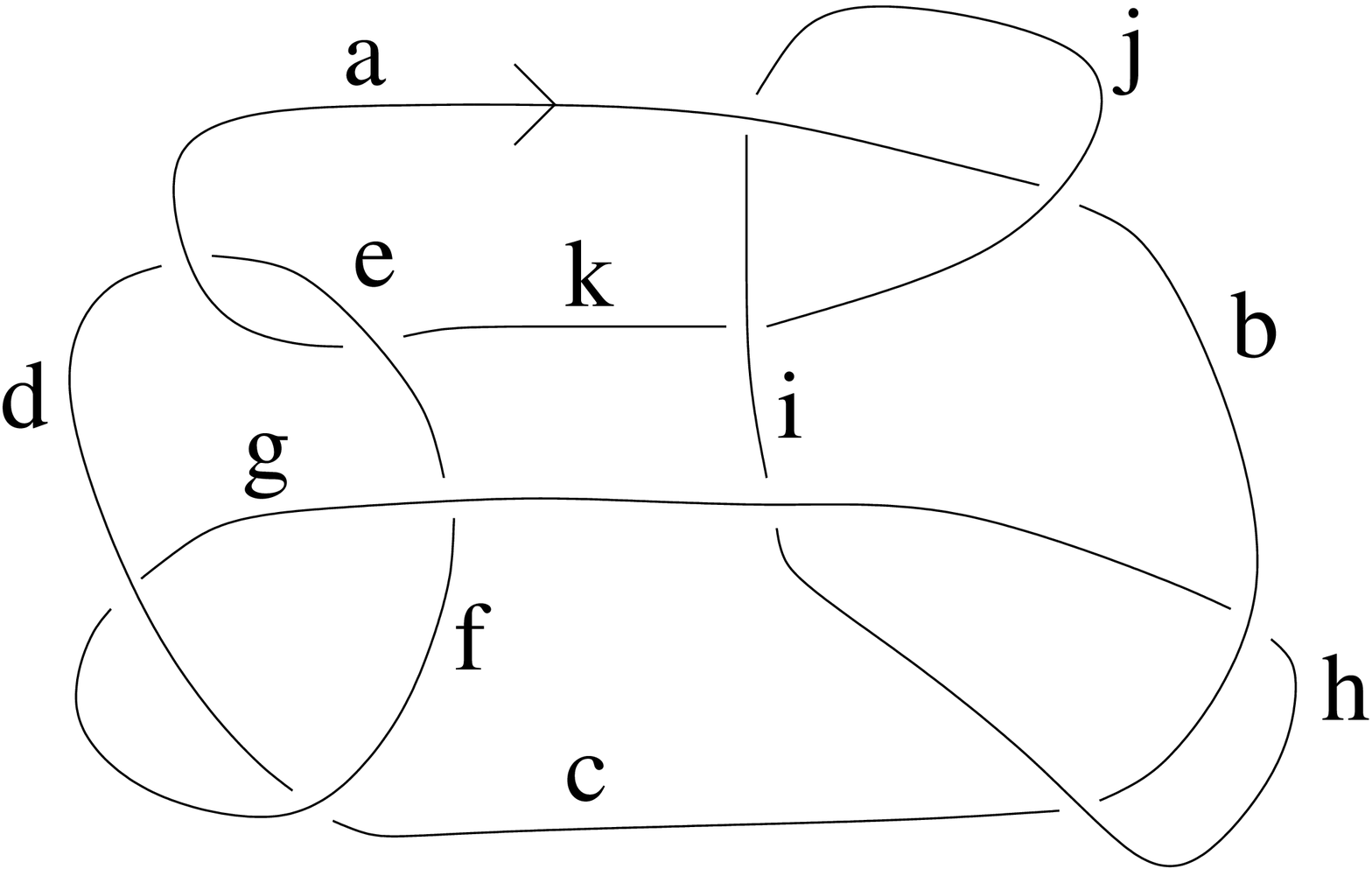} \end{tabular} \caption{Link $L$ and knot
$K_2=11_{440}$ with meridians.} \label{link11n73}
 \end{center}
 \end{figure}

Now consider the Hopf-like link $L$ in Figure \ref{link11n73}.
This consists of the knot $K_1$, the trefoil, and $K_2=11_{440}$
(here we use the \emph{knotscape} notation). By Corollary
\ref{coralexlink} the usual multivariable Alexander polynomial
with rational coefficients equals
\[ \Delta_L(x_1,x_2)=\Delta_{K_1}(x_1)\Delta_{K_2}(x_2)=(x_1^2-x_1+1)(x_2^4-2x_2^3+3x_2^2-2x_2+1)\in \Q
\xypm.\]
Let $\{\phi_1,\phi_2\} \subset H^1(X(L);\Z) = \hom(H_1(X(L);\Z),\Z)$
be the dual basis to $\{x_1,x_2\}$. It is known that $\genus(K_1)=1$
and $\genus(K_2)=3$.  We can arrange the minimal Seifert surfaces
such that they are punctured once by the other component. It follows
that $||\phi_1||_T \leq 2\, \genus(K_1)=2$ and $||\phi_2||_T \leq
2\, \genus(K_2) = 6$. In fact it is easy to see that the equality
holds for each case since each surface dual to $\phi_1$
(respectively $\phi_2$) becomes a Seifert surface for $K_1$
(respectively $K_2$) after adding one or more disks. On the other
hand it follows from the calculation of $\Delta_L(x_1,x_2)$ that
$||\phi_1||_A=2$ and $||\phi_2||_A=4$. Therefore the Alexander norm
and the Thurston norm do not agree for $L$. We also note that since
$H_1(X(L);\Z)$ is torsion-free, Turaev's torsion norm \cite{Tu02a}
agrees with the Alexander norm.

The fundamental group of $\pi_1(X(K_2))$ is generated by the
meridians $a,b,\dots,k$ of the segments in the knot diagram in
Figure \ref{link11n73}. Using the program \emph{KnotTwister}
\cite{F05} we found the homomorphism $\varphi:\pi_1(X(K_2))\to
S_3$ given by
 \[
\ba{rclrclrclrclrclrcl}
 A&=&(23),& B&=&(12),& C&=&(13),& D&=&(23), &E&=&(23),& F&=&(12),\\
 G&=&(13),& H&=&(23),& I&=&(12),& J&=&(13),& K&=&(23),&\ea \]
 where we use the cycle notation. The generators of $\pi_1(X(K_2))$ are sent to the elements in
$S_3$ given by the cycle with the corresponding capital letter. We then consider
$\a:=\a(\varphi):\pi_1(X(K_2))\xrightarrow{\varphi} S_3\to \gl(V_2)$ where
\[
V_2:=\{(v_1,v_2,v_3)\in \F_{13}^{3} | \sum_{i=1}^{3} v_i =0 \}.
\]
Clearly $\dim_{\f_{13}}(V_2)=2$ and $S_{3}$ acts on it by
permutation. With \emph{KnotTwister} we compute
\[ \Delta^\a_{K_2}(x_2)=1+3x_2^2+12x_2^4+x_2^6+10x_2^8+12x_2^{10} \in \F_{13}\ypm \]
and $H_0^\a\left(X(K_2);\F_{13}^2[x_2^{\pm 1}]\right)=0$. Hence
$\Delta^{\a,0}_{K_2}(x_2)=1$.

\begin{figure}[h] \begin{center}
\begin{tabular}{rl}
\includegraphics[scale=0.3]{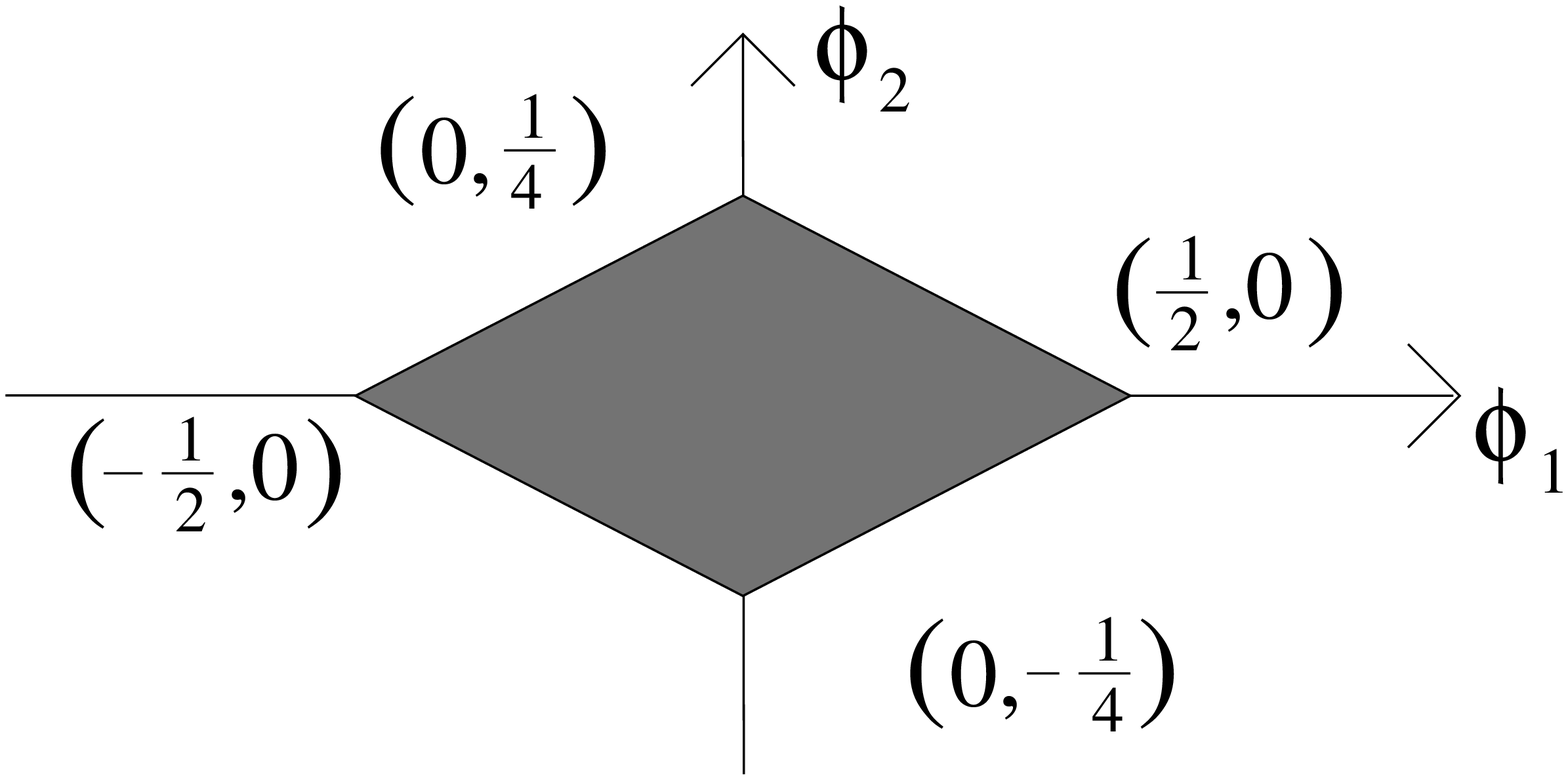}&\includegraphics[scale=0.3]{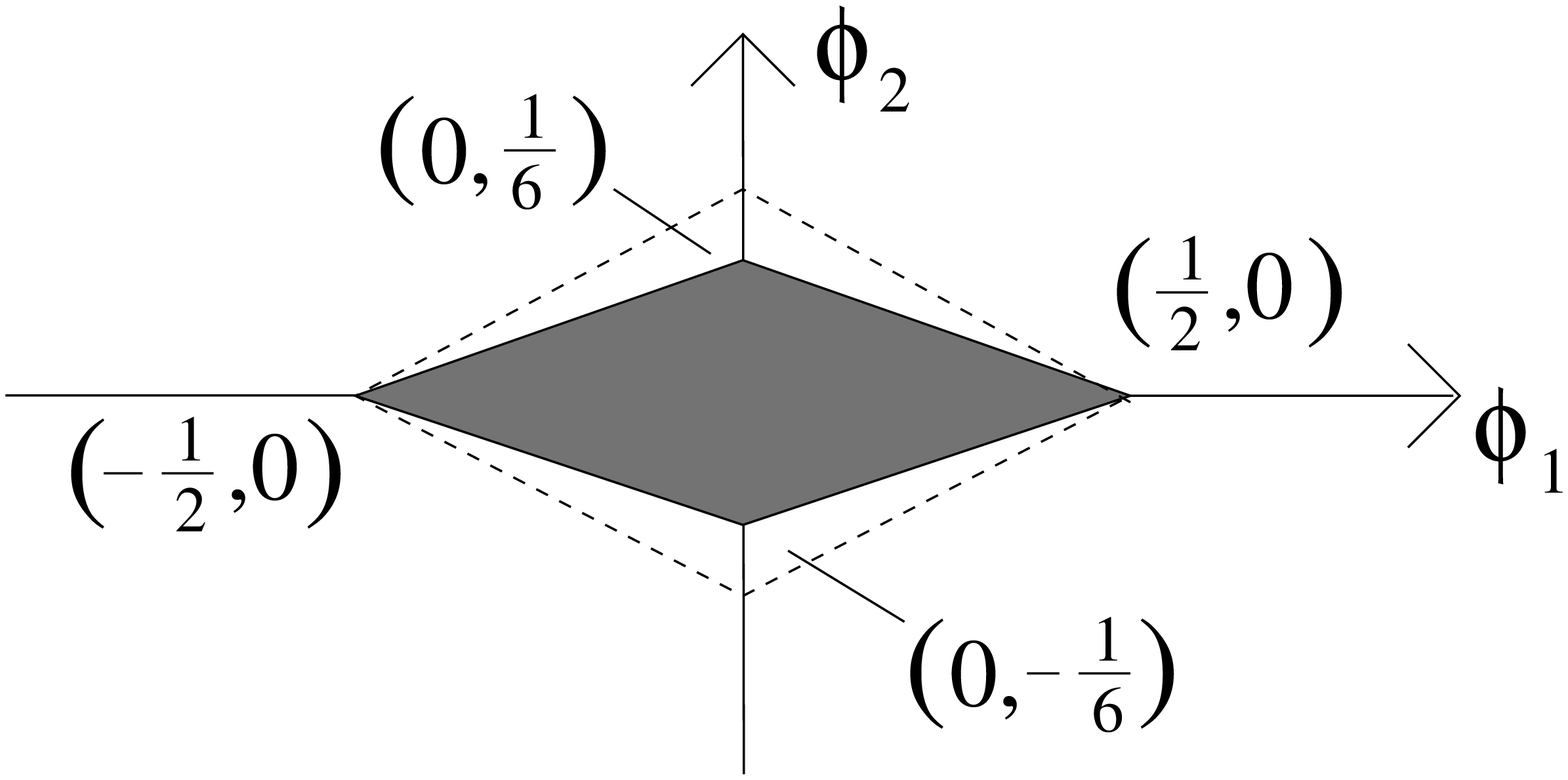}
\end{tabular}
 \caption{The untwisted and the twisted Alexander norm ball of $L$.} \label{normball}
\end{center} \end{figure}

 Denote the homomorphism $\a:\pi_1(X(L))\to \pi_1(X(K_2))\to \gl(V_2)$ by $\a$ as
well. Here the map $\pi_1(X(L)) \to \pi_1(X(K_2))$ is induced from the inclusion. This induces a
representation of $\pi_1(X(K_1))$ as in the proof of Proposition \ref{propalexlink}, and we also
denote it by $\a$. In fact, one easily sees that $\a : \pi_1(X(K_1)) \to \gl(V_2)$ is trivial. This
implies that $\Delta^\a_{K_1}(x_1) = (\Delta_{K_1}(x_1))^2 = (1-x_1 + x_1^2)^2$ and
$\Delta^{\a,0}_{K_1}(x_1)= (x_1-1)^2$. By Proposition \ref{propalexlink} we have
\[
\Delta_L^\a(x_1,x_2)=\Delta_{1}^\a(x_1)\cdot \Delta^\a_{2}(x_2)
\]
where
\[
\deg(\Delta_1^\a(x_1))=2\, \deg\left(\Delta_{K_1}(x_1)\right)+2-2
= 4
\]
and
\[
\deg(\Delta_2^\a(x_2))=\deg\left(\Delta^\a_{K_2}(x_2)\right)+2-0 =
12.
\]
Hence the twisted Alexander norm ball corresponding to
$\frac{1}{2}||-||_A^\a$ has exactly four extreme vertices
$(\pm\frac12,0)$ and $(0,\pm\frac16)$ by Corollary
\ref{corvertex}. Since $||\phi_1||_T = 2$ and $||\phi_2||_T = 6$,
the norms $||\phi||_T$ and $\frac{1}{2} ||\phi||_{A}^\a$ agree at
the extreme vertices of the norm ball of $\frac12||-||_{A}^\a$.
Note that by Theorem \ref{mainthm} we have $||\phi||_T \ge
\frac{1}{2} ||\phi||_{A}^\a$. Since the norms $||\phi||_T$ and
$\frac{1}{2} ||\phi||_{A}^\a$ agree at all of the extreme vertices
of the norm ball of $\frac12 ||-||_{A}^\a$, they agree everywhere
by convexity. Therefore the shaded region on the right in Figure
\ref{normball} is the Thurston norm ball of the link $L$.

In Figure \ref{normball} on the right the closed region bounded by the dashed polygon is the
Alexander norm ball. If $(X(L),\phi)$ fibers over $S^1$ for some $\phi\in H^1(X(L);\Z)$ then it
follows from Theorem \ref{mainthmfib} that the (usual) Alexander norm and the Thurston norm agree
on the cone on a top-dimensional face of the Thurston norm ball. Figure \ref{normball} shows that
the Alexander norm and the Thurston norm agree only for a multiple of $\phi_1$. Hence $(X(L),\phi)$
does not fiber over $S^1$ for any $\phi \in H^1(X(L);\Z)$. We state these results   in the
proposition below.
\begin{proposition} \label{thmhophlikelink}
The Thurston norm ball of $X(L)$ is the shaded region on the right
in Figure \ref{normball}. Furthermore, $(X(L), \phi)$ does not
fiber over $S^1$ for any $\phi \in H^1(X(L);\Z)$.
\end{proposition}


There exist 36 knots with 12 crossings or less such that $2\,
\genus(K)>\deg(\Delta_K(t))$. In all but three cases we found
representations similar to the above such that the Thurston norm
bound from Theorem \ref{mainthmfk05} equals the Thurston norm of
$X(K)$. Let $L$ be the Hopf-like link as in Figure \ref{hopflink}
with  $K_1$ any knot such that $2\,
\genus(K_1)=\deg(\Delta_{K_1}(t))$ and $K_2$ any of the 33 knots
mentioned above. In this case the argument above can be used to show
that twisted Alexander norms completely determine the Thurston norm
ball of $X(L)$ and it is always strictly smaller than the Alexander
norm ball.
\\

\begin{figure}[h]
\begin{center}
\includegraphics[scale=0.3]{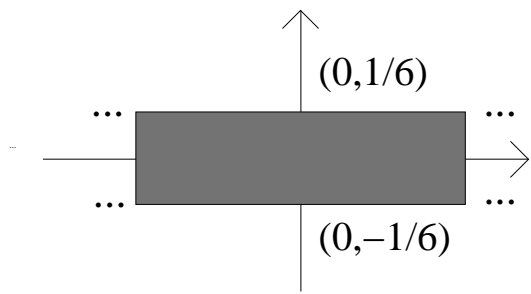}
\caption{Thurston norm ball of $L$.} \label{normball2}
\end{center}
\end{figure}

Now consider the case with $K_1$ the unknot and $K_2=11_{440}$. We
use the same representation as above. In this case the norm ball
for $\frac{1}{2}||-||_A^\a$ is given in Figure \ref{normball2}.
The norm ball is a horizontal infinite strip, hence noncompact.

To show that $\frac{1}{2}||-||_A^\a=||-||_T$ it is enough to show
that for $\phi=(n,\pm 1), n\in \Z$ there exists a connected dual
surface with $\chi(S)=-6$. Let $S$ be a Seifert surface of genus 3
for $K_2$ which intersects $K_1$ just once. By deleting a disk
from $S$ we get a surface $S'$ which is disjoint from $K_1$. The
surface $S'$ is dual to $\phi = (0,1)$. We can make $S'$ such that
the two boundary components of $S'$ are as close to each other as
we wish. Now take a short path from one boundary component of $S'$
to the other boundary component. Cut $S'$ along that path and
reglue the cut parts together by giving $n$ full twists. The
resulting surface is dual to $\phi=(n,1)$ and has the Euler
characteristic -6. Hence the Thurston norm ball in this case is
the shaded (infinite) strip in Figure \ref{normball2}.

\subsection{Dunfield's example} \label{sec:dunfield}
\begin{figure}[h] \begin{center} \includegraphics[scale=0.3]{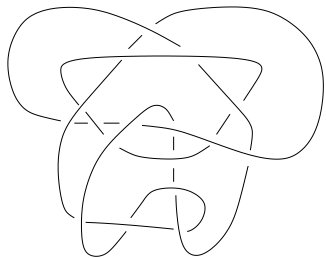}
\caption{Dunfield's example.} \label{dunfield} \end{center}
\end{figure}

McMullen had asked whether for a fibered manifold the Thurston norm and the Alexander norm agree
everywhere. To answer this question Dunfield \cite{Du01} considers the link $L$ in Figure
\ref{dunfield}.

\begin{figure}[h] \begin{center}
\includegraphics[scale=0.3]{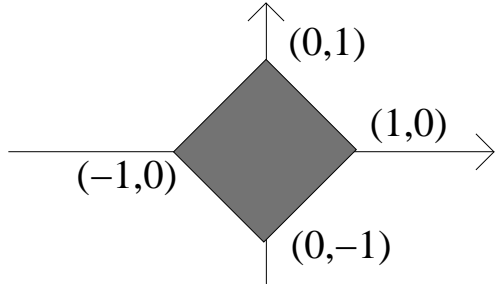} \caption{Alexander norm ball for Dunfield's link.}
\label{dunfieldanorm} \end{center}
 \end{figure}

Denote the knotted component by $K_1$ and the unknotted component
by $K_2$. Let $x,y\in H_1(X(L);\Z)$ be the elements represented by
a meridian of $K_1$, respectively  $K_2$. Then the Alexander
polynomial equals
 \[ \Delta_{X(L)}=xy-x-y+1 \in \Z[H_1(X(L);\Z)]=\Z[x^{\pm 1},y^{\pm 1}].\]
We consider $H^1(X(L);\Z)$ with the dual basis corresponding to
$\{x,y\}\in H_1(X(L);\Z)$. The Alexander norm ball is given in
Figure \ref{dunfieldanorm}.

Dunfield \cite{Du01} showed  that $(X(L),\phi)$ fibers over $S^1$
for all $\phi\in H^1(X(L);\Z)$ in the cones on the two open faces
of the Alexander norm ball with vertices
$(-\frac{1}{2},\frac{1}{2}), (0,1)$ respectively $(0,-1),
(\frac{1}{2},-\frac{1}{2})$. Dunfield used the
Bieri-Neumann-Strebel (BNS) invariant (see \cite{BNS87}) to show
that the Alexander norm and the Thurston norm do not agree for the
3-manifold $X(L)$. We will go one step further and completely
determine the Thurston norm of $X(L)$.

We did not find a representation of $\pi_1(X(L))$ for which we can compute the twisted Alexander
polynomial and which determines the Thurston norm. Therefore we study the Thurston norm of a
2-fold cover of $X(L)$ for which it is easier to find a representation.

The following theorem by Gabai shows the relationship between the Thurston norm of $X(L)$ and that
of a finite cover of $X(L)$.

\begin{theorem} \label{lemmathurstong} \cite[p.~484]{Ga83}
Let $M$ be a 3-manifold and $\a:\pi_1(M)\to G$ a homomorphism to a
finite group $G$. Denote the induced $G$-cover of $M$ by $M_G$.
Let $\phi\in H^1(M;\Z)$ be nontrivial and denote the induced map
$H_1(M_G;\Z)\to H_1(M;\Z)\to \Z$ by $\phi_G$, which can be
regarded as an element in $H^1(M_G;\Z)$. Then $\phi_G$ is
nontrivial and
\[ |G|\cdot||\phi||_{T,M}= ||\phi_G||_{T,M_G}.\]
\end{theorem}
\noindent Thus to determine the Thurston norm of $M$, we only need
to determine the Thurston norm of $M_G$. For this purpose, we
generalize twisted Alexander norms and the main theorems a little
bit further as follows.


Let $M$ be a 3-manifold and $\psi:\pi_1(M)\to F$ a homomorphism to
a free abelian group, we do not demand that $\psi$ is surjective. We define a norm on $\Hom (F,\R)$. Note that
if $F = FH_1(M;\Z)$, then $\Hom (F,\R) \cong H^1(M;\R)$. Let
$\a:\pi_1(M)\to \gl (\F,k)$ be a representation. If
$\Delta_{M,\psi}^{\a}=0\in \F[F]$ then we set
$||\phi||_{A,\psi}^{\a}=0$ for all $\phi\in \Hom(F,\R)$. Otherwise
we write $\Delta_{M,\psi}^{\a}=\sum a_if_i$ for $a_i\in \F$ and
$f_i \in F$. Given $\phi \in \Hom(F,\R)$, we define \emph{the
(generalized) twisted Alexander norm of $(M,\psi,\a)$} to be
\[ ||\phi||_{A,\psi}^{\a} :=\sup \phi(f_i-f_j)\]
with the supremum over $(f_i, f_j)$ such that $a_ia_j\ne 0$. If we consider the natural surjection
$\psi:\pi_1(M)\to FH_1(M;\Z)$, then clearly $||-||_{A,\psi}^{\a}=||-||_A^{\a}$. Note that
$||-||_{A,\psi}^{\a}$ is clearly a seminorm on $\hom(F,\R)$. The following theorem generalizes
Theorem \ref{mainthm} and Theorem \ref{mainthmfib}. The proof is almost identical.

\begin{theorem} \label{thmgeneral}
Let $M$ be a 3-manifold whose boundary is empty or consists of tori.
Let $\a:\pi_1(M)\to \glfk$ be a representation. Let
$\psi:\pi_1(M)\to F$ be a homomorphism to a free abelian group such
that $\rank\, F>1$ and such that $H_1(M;\Z)\otimes_\Z \Q \to
F\otimes_\Z \Q$ is surjective. Then
\[ ||\phi\circ \psi ||_T \ge \frac1k ||\phi||^\a_{A,\psi}   \]
for all $\phi\in \Hom(F,\R)$.

Furthermore, if $M\ne S^1\times D^2, M \ne S^1\times S^2$ and if $\phi\in \Hom(F,\Z)$ is such that
$(M,\phi\circ \psi)$ fibers over $S^1$, then $\phi\circ \psi$ lies in the cone on a top-dimensional
open face of the Thurston norm ball (denoted by $C$) and for all $\phi'\in \Hom(F,\R)$ such that
$\phi'\circ \psi \in C$ we have
\[ ||\phi'\circ \psi||_T = \frac1k ||\phi'||^\a_{A,\psi}.  \]
\end{theorem}


We now return to the link $L$ in Figure \ref{dunfield}. Let
$\varphi:H_1(X(L);\Z)\to \Z/2$ be the homomorphism given by
$\varphi(x)=1$, $\varphi(y)=0$. Denote the induced two-fold cover by
$X(L)_2$. Denote by $\psi$ the homomorphism $\pi_1(X(L)_2) \to
H_1(X(L)_2;\Z)\to H_1(X(L);\Z)$ induced from the covering map $\pi :
X(L)_2 \to X(L)$. We found a representation $\a:\pi_1(X(L)_2)\to
\gl(\F_7,1)$ such that
\[ \Delta^\a_{X(L)_2,\psi} = 3x^6y^2+3x^4y^2+4x^4y+2x^4+x^2y^2+3x^2y-x^2-1
\in \F_7[H_1(X(L);\Z)]=\F_7[x^{\pm 1},y^{\pm 1}].\] This polynomial is not of the form $f(ax+by)$
for some polynomial $f(t)$. This shows that $H_1(X(L_2))\to H_1(X(L))=\Z^2$ is rationally
surjective, in particular we can apply Theorem \ref{thmgeneral}.


%

\begin{figure}[h]
\begin{center}
\includegraphics[scale=0.3]{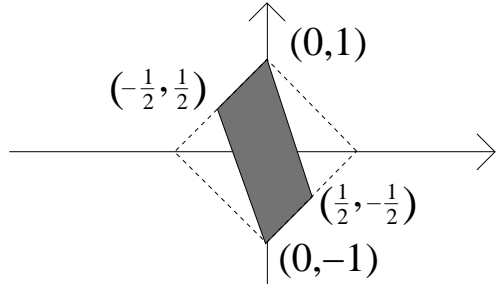}
\caption{Twisted Alexander norm ball for Dunfield's link}
\label{dunfieldtwinorm}
\end{center}
\end{figure}

Now let $\phi \in H^1(X(L);\Z)$. By Theorem \ref{lemmathurstong} and Theorem \ref{thmgeneral},  we
have
\[
||\phi||_{T,X(L)}=\frac{1}{2}||\phi \circ\pi ||_{T,X(L)_2} \ge \frac{1}{2} ||\phi||^\a_{A,\psi}.
\]
The norm ball of $\frac{1}{2} ||-||^\a_{A,\psi}$ is drawn as the shaded region in Figure
\ref{dunfieldtwinorm}. We claim that this is exactly the Thurston norm ball.

By Theorem \ref{thmgeneral} the twisted Alexander norm ball in
Figure \ref{dunfieldtwinorm} is an `outer bound' for the Thurston
norm ball of $X(L)$. But as we pointed out above, Dunfield showed
that $(X(L),\phi)$ fibers over $S^1$ for all $\phi\in
H^1(X(L);\Z)$ which lie in the cones on the two open faces of the
Alexander norm ball with vertices $(-\frac{1}{2},\frac{1}{2}),
(0,1)$ respectively $(0,-1), (\frac{1}{2},-\frac{1}{2})$. In
particular, the Thurston norm ball and the twisted Alexander norm
ball agree on these cones by the second part of Theorem
\ref{thmgeneral}. By continuity, the norms also agree on the
vertices $(-\frac{1}{2},\frac{1}{2}), (0,1), (0,-1)$ and
$(\frac{1}{2},-\frac{1}{2})$. Now it follows from convexity that
the Thurston norm ball coincides everywhere with the twisted
Alexander norm ball given in Figure \ref{dunfieldtwinorm}.
Therefore the shaded region in Figure \ref{dunfieldtwinorm} is the
Thurston norm ball of $X(L)$.

Note that our calculation confirms Dunfield's result that $(X(L),\phi)$ does not fiber over $S^1$
for any $\phi$ outside
the cones. We summarize these results in the following proposition.

\begin{proposition} \label{prop:dunfield}
The Thurston norm ball of $X(L)$ is the shaded region in Figure
\ref{dunfieldtwinorm}. Furthermore, $(X(L),\phi)$ fibers over
$S^1$ exactly when $\phi$ lies inside the cones on the open faces
of the two smaller faces of the Thurston norm ball of $X(L)$.
\end{proposition}

\section{Twisted multivariable Alexander polynomial and twisted one-variable Alexander polynomial}
\label{sec:multi}

This section serves for proving Theorem \ref{mainthmalex}. The
main idea of the proof is to use the functoriality of Reidemeister
torsion. To prove Theorem \ref{mainthmalex} we need some lemmas
which show the nontriviality of certain twisted Alexander
polynomials. Throughout this section we assume that $M$ is a
3-manifold whose boundary is empty or consists of tori.
Furthermore let $\a : \pi_1(M) \to \glfk$ be a representation.

\subsection{Computation of twisted Alexander polynomials} \label{sec:calc}
We introduce the notion of \emph{rank over a UFD}. Let $\L$ be a UFD
and $Q(\L)$ its quotient field. Let $H$ be a $\L$-module. Then we
define $\rank_\L(H):=\dim_{Q(\L)}(H\otimes_{\L}Q(\L))$. We need the
following well-known lemma. For the first part we refer to
\cite[Remark~4.5]{Tu01}. The second part is well--known. The last
statement follows from the fact that $Q(\L)$ is flat over $\L$.

\begin{lemma} \label{lemma:rank}
Let $\L$ be a UFD. \bn \item Let $H$ be a finitely generated
$\L$-module. Then the following are equivalent: \bn \item $H$ is
$\L$-torsion, \item $\ord_{\L}(H)\ne 0$, \item $\rank_\L(H)=0$,
\item $\hom_{\L}(H,\L)=0$. \en \item Let $N$ be an n-manifold and
assume that $\L^k$ has a left $\Z[\pi_1(N)]$-module structure,
then
\[ \sum \limits_{i=0}^n (-1)^i \rank_\L(H_i(N;\L^k))=k\chi(N).\]
\item $H_i(N;\L^k\otimes_\L {Q(\L)})=H_i(N;\L^k)\otimes_\L
{Q(\L)}$ for any $i$. \en
\end{lemma}

\begin{lemma} \label{lem:delta03}
Let $M$ be a 3-manifold. Let $\varphi : \pi_1(M) \to H$ be a
surjection to a free abelian group.  Then
$\Delta_{M,\varphi}^{\a,3}=1$ and $\Delta_{M,\varphi}^{\a,0}\ne
0$. If furthermore  $\rank\, H > 1$, then $
\Delta_{M,\varphi}^{\a,0}=1 \in \F[H]$.
\end{lemma}

\begin{proof}
We prove the lemma only in the case that $M$ is closed. The proof for the case that $\partial M$
consists of tori is very similar. Let $b:=\rank\, H$. Pick a basis $t_1,\dots,t_b$ for $H$.  We
identify $\F^k[H]:=\F^k \otimes \F[H]$ with $\fk[t_1^{\pm 1},\dots,t_b^{\pm 1}]$.

Since $M$ is closed it follows that $\chi(M)=0$. Then it is
well--known that $M$ has a CW--structure with one cell of dimensions
zero and three, and with the same number of cells in dimensions one
and two (cf. e.g. \cite[Theorem 5.1]{Mc02}).
Denote the 1-cells by $h_1,\dots,h_n$. Denote the corresponding elements in $\pi_1(M)$ by
$h_1,\dots,h_n$ as well. If $\rank\, H>1$ then we can arrange that $\varphi(h_i)=t_i$ for $i=1,2$.

Write $\pi:=\pi_1(M)$. From the CW structure we obtain a chain
complex $C_*:=C_*(\ti{M})$ (where $\ti{M}$ denotes the universal
cover of $M$):
\[
0 \to C_3^1 \xrightarrow{\partial_3} C_2^n \xrightarrow{\partial_2} C_1^n \xrightarrow{\partial_1}
C_0^1 \to 0
\]
for $M$, where the $C_i$ are free $\Z[\pi]$-right modules. In fact $C_i^k\cong \Z[\pi]^k$. Consider
the chain complex $C_*\otimes_{\Z[\pi]}\F^k[H]$:
\[
0 \to C_3^1\otimes_{\Z[\pi]}\fk[H] \xrightarrow{\partial_3\otimes
\id} C_2^n\otimes_{\Z[\pi]}\fk[H] \xrightarrow{\partial_2\otimes
\id} C_1^n \otimes_{\Z[\pi]}\fk[H] \xrightarrow{\partial_1\otimes
\id} C_0^1\otimes_{\Z[\pi]}\fk[H]  \to 0.
\]

Let $A_i$, $i=0,\dots,3$, be the matrices with entries in $\Z[\pi]$ corresponding to the boundary
maps $\partial_i:C_i\to C_{i-1}$ with respect to the bases given by the lifts of the cells of $M$
to $\ti{M}$. Then $A_3$ and $A_1$ are well--known to be of the form
 \[ \ba{rcl} A_3 &=&
(a_1(1-g_1), a_2(1-g_2), \ldots, a_n(1-g_n))^t,\\
A_1 &=& (b_1(1-h_1), b_2(1-h_2), \ldots, b_n(1-h_n)), \ea \] where
$\{g_1,\dots,g_n\}$ and $\{h_1,\dots,h_n\}$ are generating sets for
$\pi_1(M)$ and $a_i,b_i \in \pi_1(M)$ for $i=1,\dots,n$. By picking
different lifts of the cells in dimensions one and two we can assume
that in fact $a_i=b_i=e\in \pi_1(M)$ for $i=1,\dots,n$. We can and
will therefore assume that
 \[ \ba{rcl} A_3 &=&
(1-g_1, 1-g_2, \ldots, 1-g_n)^t,\\
A_1 &=& (1-h_1, 1-h_2, \ldots, 1-h_n). \ea \]

Let $B = (b_{rs})$ be a $p\times q$ matrix with entries in
$\Z[\pi]$. We write $b_{rs}=\sum b_{rs}^gg$ for $b_{rs}^g\in \Z,
g\in \pi$. We define $(\a\otimes \varphi)(B)$ to be the $p\times q$
matrix with entries $\sum b_{rs}^g \a(g)\varphi(g)$. Since each
$\sum b_{rs}^g \a(g)\varphi(g)$ is a $k\times k$ matrix with
entries in $\F[H]$ we can think of  $(\a\otimes \varphi)(B)$ as a
$pk\times qk$ matrix with entries in $\F[H]$.

Since $\varphi$ is nontrivial there exist $k,l$ such that
$\varphi(g_k)\ne 0$ and $\varphi(h_l)\ne 0$. It follows that $(\a
\otimes \varphi)(A_1)$ and $(\a \otimes \varphi)(A_3)$ have full
rank over $\f[H]$. The first part of the lemma now follows
immediately.

Now assume that $\rank\,H>1$. Then $\ord\left(H_0^\a(M;\fk[t_1^{\pm
1},\dots,t_b^{\pm 1}])\right)$ divides $\det(\a(h_1)t_1-\id)\in \f
[t_1^{\pm 1}]$ and $\det(\a(h_2)t_2-\id)\in \f[t_2^{\pm 1}]$. These
two polynomials are clearly  relatively prime. This implies that
$\ord\left(H_0^\a(M;\fk[t_1^{\pm 1},\dots,t_b^{\pm 1}])\right)=1$.

\end{proof}

\begin{lemma} \label{lem:delta2}
Let $\varphi:\pi_1(M)\to H$ be a surjection to a free abelian group $H$.
If $\Delta_{M,\varphi}^{\a,1}\ne 0$ then
$\Delta_{M,\varphi}^{\a,2}\ne 0$.
\end{lemma}

\begin{proof} Note that by assumption and by Lemma \ref{lem:delta03} we have $\Delta_{M,\vph}^{\a,i}
\ne 0$ for $i=0,1,3$. Let $\L:=\F[H]$. It follows from the long
exact homology sequence for $(M,\partial M)$ and from duality that
$\chi(M)=\frac{1}{2}\chi(\partial M)$. So $\chi(M)=0$ in our case.
It follows from Lemma \ref{lemma:rank} that
\[ \sum_{i=0}^3 (-1)^i \dim_{Q(\L)} \left( H_i^\a(M;\L^k \otimes_\L Q(\L)   )\right)=k\chi(M)=0.\]
Note that $H_i^\a(M;\L^k\otimes_\L Q(\L)) \cong
H_i^\a(M;\L^k)\otimes_\L Q(\L)$ by Lemma \ref{lemma:rank}. By
assumption $H_i^\a(M;\L^k)\otimes_{\L} Q(\L)=0$ for $i\ne 2$,
hence $H_2^\a(M;\L^k)\otimes_\L Q(\L)=0$.
\end{proof}

The following corollary is now immediate.

\begin{corollary} \label{coracyclic}
Let $\varphi:\pi_1(M)\to H$ be a surjection to a free abelian
group $H$. If $\Delta_{M,\varphi}^{\a,1}\ne 0$ then
$\Delta_{M,\varphi}^{\a,i}\ne 0$ for all $i$.
\end{corollary}

\begin{lemma} \label{lem:delta2multi}
Let $\varphi : \pi_1(M) \to H$ be a surjection to a free abelian
group with $\rank\, H > 1$. If $\Delta_{M,\varphi}^{\a,1} \ne 0$
then
\[ \Delta_{M,\varphi}^{\a,2} =1 \in \F[H].\]
\end{lemma}

\begin{proof}
Let $\L:=\F[H]$ and $\pi :=\pi_1(M)$.
 By Poincar\'e
duality,
\[ H_2^\a(M;\L^k) \cong H^1_\a(M,\partial M;\L^k)=H^1(\Hom_{\Z[\pi]}(C_*(\tilde{M},\partial \ti{M}), \L^k))
\]
where $\ti{M}$ is the universal cover of $M$. On the right we view
$\L^k$ as a right $\Z[\pi]$-module by taking $f\cdot g =
g^{-1}\cdot f=\varphi(g^{-1})\a(g^{-1})f$ for $f\in \L^k$ and $g\in \pi$.

We use an argument in \cite[p.~638]{KL99}. Let $\langle \, , \,
\rangle : \fk \times \fk \to \f$ be the canonical inner product on
$\fk$.
Then
there exists a unique representation $\overline{\a}:\pi_1(M)\to
\glfk$ such that
\[ \langle \a(g^{-1})v,w \rangle = \langle v,\overline{\a}(g)w\rangle\]
for all $g\in \pi_1(M)$ and $v,w\in \fk$. We denote by
$\overline{\L^k}$ the left $\Z[\pi]$-module with underlying
$\L$-module $\L^k$ and $\Z[\pi]$-module structure given by $\at
\otimes (-\phi)$.

Using the inner product we get a map
\[\ba{rcl}\Hom_{\Z[\pi]}(C_*(\tilde{M},\partial \ti{M}), \L^k)&\to&
\hom_{\L}\big(C_*(\ti{M},\partial \ti{M})\otimes_{\Z[\pi]} \overline{\L^k},\L\big)\\
f&\mapsto& (c\otimes w)\mapsto \langle f(c),w\rangle. \ea \] Using
$\langle \a(g^{-1})v,w\rangle =\langle v,\at(g)w\rangle$ it is now
easy to see that this map is well-defined and that it defines in
fact an isomorphism of $\L$-module chain complexes.

Now we can apply the universal coefficient spectral sequence to
the $\L$-module chain complex $\hom_{\L}\big(C_*(\ti{M},\partial
\ti{M})\otimes_{\Z[\pi]} \overline{\L^k},\L\big)$ to conclude that
there exists a short exact sequence
\[ 0\to \ext_{\L}^1(H_0^{\at}(M,\partial M;\overline{\L^k}))\to H^1_\a(M,\partial M;\L^k)
\to  \Hom_{\L}(H_1^{\at}(M,\partial M;\overline{\L^k})).\]
 Since $\Delta^{\a,2}_{M,\varphi}\ne 0$ by Lemma \ref{lem:delta2}
it follows that $H^1_\a(M,\partial M;\L^k)$ is $\L$-torsion. Hence
\[ H^1_\a(M,\partial M;\L^k)  \cong \ext_{\L}^1(H_0^{\at}(M,\partial M;\overline{\L^k}).\]
First assume that $\partial M$ is nonempty. Note that
$\pi_1(\partial M)\to \glfk$ factors through $\pi_1(M)$. It
follows from
\[ H_0^{\at}(X;\overline{\L^k})\cong \overline{\L^k} /\{ gv-v | g\in \pi_1(X), v\in \overline{\L^k}\}\]
that ${H}_0^{\at}(\partial M;\overline{\L^k})$ surjects onto ${H}_0^{\at}(M;\overline{\L^k})=0$,
hence ${H}_0^{\at}(M,\partial M;\overline{\L^k})=0$ (cf. \cite[Lemma~2.6]{FK05}).

Now assume that $M$ is closed. Let $H_0 := {H}_0^{\at}(M;\overline{\L^k})$. We define a finitely
generated $\L$-module $A$ to be \emph{pseudonull} if $A_\wp = 0$ for every height 1 prime ideal
$\wp$ of $\L$ where $A_\wp$ is the localization of $A$ at $\wp$. (See p.~51 in \cite{Hi02}.) By
\cite[Theorem 3.1]{Hi02}, $E_0(H_0) \subset \ann(H_0)$. Since $\Delta^{\a,0}_{M,\vph} = 1$ by Lemma
\ref{lem:delta03}, $\widetilde{\ann}(H_0) = \L$ where $\widetilde{\ann}(H_0)$ is the smallest
principal ideal of $\L$ which contains $\ann(H_0)$. Thus by \cite[Theorem 3.5]{Hi02}, $H_0$ is
pseudonull. Finally, by \cite[Theorem 3.9]{Hi02}, $\Ext^1_\L(H_0, \L) = 0$. Hence $H_2^\a(M;\L^k)
\cong H^1_\a(M,\partial M;\L^k)= 0$.
\end{proof}

\subsection{Functoriality of torsion} \label{sec:funtoriality}
Define $F$ to be the free abelian group $FH_1(M;\Z)$. Let $\psi :
\pi_1(M) \to F$ be the natural surjection and $\phi \in H^1(M;\Z)$
nontrivial. Note that $\phi$ induces a homomorphism $\phi : \F[F]
\to \F[t^{\pm 1}]$. In this section we go back to the notation
$\Delta_M^{\a,i} = \Delta_{M,\psi}^{\a,i}$ and $\Delta_\phi^{\a,i}
= \Delta_{M,\phi}^{\a,i}$.

\begin{theorem} \label{thm:turaev}
Suppose $b_1(M) > 1$. \bn  \item If
$\phi\left(\Delta_{M,\psi}^{\a,1}\right) \ne 0$ then
$\Delta_{M,\phi}^{\a,1} \ne 0$ and
\[ \phi\left(\Delta_{M,\psi}^{\a,1}\right) =  \prod \limits_{i=0}^3 \phi\left(\Delta_{M,\psi}^{\a,i}\right)^{(-1)^{i+1}}=\prod \limits_{i=0}^3
\Delta_{M,\phi}^{\a,i}(t)^{(-1)^{i+1}}\in \ft.
\]
\item If $\phi\left(\Delta_{M,\psi}^{\a,1}\right) = 0$ then
$\Delta_{M,\phi}^{\a,1} = 0$. \en
\end{theorem}

\noindent Note that if $\Delta_{M,\psi}^{\a,1} \ne 0$, then by
Lemmas \ref{lem:delta03} and \ref{lem:delta2multi} $\prod_{i=0}^3
\phi\left(\Delta_{M,\psi}^{\a,i}\right)^{(-1)^{i+1}}$ is defined and
the first equality in the first part is obvious. Also if
$\Delta_{M,\phi}^{\a,1} \ne 0$ then by Lemmas \ref{lem:delta03} and
\ref{lem:delta2}, $\prod_{i=0}^3
\Delta_{M,\phi}^{\a,i}(t)^{(-1)^{i+1}}$ is defined.

\begin{proof}
We will only consider the case that $M$ is a closed 3-manifold. The proof for the case that $\partial M\ne
\emptyset$ is similar.

Let us prove (1). Write $\pi:=\pi_1(M)$. As in the proof of Lemma
\ref{lem:delta03} \label{lem:delta03} we can find a CW-structure
for $M$ such that the chain complex $C_*(\ti{M})$ of the universal
cover is of the form
\[
0 \to C_3^1 \xrightarrow{\partial_3} C_2^n
\xrightarrow{\partial_2} C_1^n \xrightarrow{\partial_1} C_0^1 \to
0
\]
for $M$, where the $C_i$ are free $\Z[\pi]$-right modules. In fact
$C_i^k\cong \Z[\pi]^k$. Let $\varphi:\pi_1(M)\to H$ be an
epimorphism to a free abelian group $H$. Consider the chain
complex $C_*\otimes_{\Z[\pi]}\F^k[H]$:
\[
0 \to C_3^1\otimes_{\Z[\pi]}\F^k[H] \xrightarrow{\partial_3\otimes
\id} C_2^n\otimes_{\Z[\pi]}\F^k[H] \xrightarrow{\partial_2\otimes
\id} C_1^n \otimes_{\Z[\pi]}\F^k[H]\xrightarrow{\partial_1\otimes
\id} C_0^1\otimes_{\Z[\pi]}\F^k[H] \to 0.
\]
Lifting the cells of $M$ to $\ti{M}$ makes $C_*$ a based complex.
Denote the quotient field of $\F[H]$ by $Q(H)$. If
 \[ C_*\otimes_{\Z[\pi]}Q(H)^k:=C_*\otimes_{\Z[\pi]}\F^k[H]
\otimes_{\F[H]} Q(H)\] is acyclic, then we can define the
Reidemeister torsion $\tau(M,\a,\varphi)\in Q(H)\sm \{0\}$ which
is well-defined up to multiplication by a unit in $\F[H]$. We
refer to \cite{Tu01} for the definition of Reidemeister torsion
and its properties.
%
%
Let $A_i$, $i=0,\dots,3$, be the matrices with entries in
$\Z[\pi]$ corresponding to the boundary maps $\partial_i:C_i\to
C_{i-1}$ with respect to the bases given by the lifts of the cells
of $M$ to $\ti{M}$. Then  we can arrange the
lifts such that
 \[ \ba{rcl} A_3 &=&
(1-g_1, 1-g_2, \ldots, 1-g_n)^t,\\
A_1 &=& (1-h_1, 1-h_2, \ldots, 1-h_n), \ea \] where
$\{g_1,\dots,g_n\}$ and $\{h_1,\dots,h_n\}$ are generating sets
for $\pi_1(M)$. Since $\phi$ is nontrivial there exist $k,l$ such
that $\phi(g_k)\ne 0$ and $\phi(h_l)\ne 0$. Let $B_3$ be the
$k$-th row of $A_3$. Let $B_2$ be the result of deleting the
$k$-th column and the $l$-th row. Let $B_1$ be the $l$-th column
of $A_1$.


Note that
\[ \det((\a\otimes \phi)(B_3))=\det(\id-(\a\otimes \phi)(g_k))
=\det(\id-\phi(g_k)\a(g_k)) \ne 0 \in \F\tpm \]
 since $\phi(g_k)\ne 0$. Similarly $\det((\a\otimes
\phi)(B_1))\ne 0$ and $\det((\a\otimes \psi)(B_i))\ne 0,  i=1,3$. We need the following theorem.
Note that $C_*\otimes_{\Z[\pi]}Q(H)$ is  acyclic if and only if $\Delta_{M,\varphi}^{\a,1}\ne 0$ by
Corollary \ref{coracyclic}.

\begin{theorem}\cite[Theorem~2.2, Lemma~2.5 and Theorem~4.7]{Tu01} \label{thm:Tu22}
Let $\varphi:\pi\to H$ be a homomorphism to a free abelian group.
Suppose $\det((\a\otimes \varphi)(B_i))\ne 0$, $i=1,3$. \bn \item
$C_*\otimes_{\Z[\pi]}Q(H)^k$ is acyclic $\Leftrightarrow
\det((\a\otimes \varphi)(B_2))\ne 0 \Leftrightarrow
\Delta_{M,\varphi}^{\a,1} \ne 0$. \item If
$C_*\otimes_{\Z[\pi]}Q(H)^k$ is acyclic then
\[ \tau(M,\a,\varphi) = \prod\limits_{i=1}^3 \det((\a\otimes \varphi)(B_i))^{(-1)^{i+1}}
=\prod\limits_{i=0}^3
\left(\Delta_{M,\varphi}^{\a,i}\right)^{(-1)^{i+1}} .\] \en
\end{theorem}

By Theorem \ref{thm:Tu22} we only need to prove that
$C_*\otimes_{\Z[\pi]}Q(F)^k$ and $C_*\otimes_{\Z[\pi]}\F(t)^k$ are
acyclic and $\tau(M,\a,\phi)=\phi(\tau(M,\a,\psi))$. (We define
$\phi(f/g) := \phi(f)/\phi(g)$ for $f,g \in \F[F]$.)

Since $\phi\left(\Delta_{M,\psi}^{\a,1}\right)\ne 0$ by our
assumption, $\Delta_{M,\psi}^{\a,1}\ne 0$. Therefore
$C_*\otimes_{\Z[\pi]}Q(F)^k$ is acyclic by Corollary
\ref{coracyclic}. Since $\det((\a\otimes \psi)(B_i))\ne 0, i=1,3$,
it follows from Theorem \ref{thm:Tu22} that $\det((\a\otimes
\psi)(B_2))\ne 0$ and
\[ \tau(M,\a,\psi)=\prod\limits_{i=1}^3 \det((\a\otimes \psi)(B_i))^{(-1)^{i+1}}.\]
Note that
\[ \ba{rcl} \prod\limits_{i=1}^3 \det((\a\otimes \phi)(B_i))^{(-1)^{i+1}}&=& \prod\limits_{i=1}^3
\phi\big(\det((\a\otimes \psi)(B_i))\big)^{(-1)^{i+1}}\\
&=& \prod\limits_{i=0}^3
\phi \left(\Delta_{M,\psi}^{\a,i}\right)^{(-1)^{i+1}}\\
&=&\phi(\tau(M,\a,\psi)). \ea
\]
In the above the second equality follows from Theorem
\ref{thm:Tu22}. Since $\phi(\Delta^{\a,1}_{M,\psi})\ne 0$ and
$\det((\a\otimes \phi)(B_i))\ne 0$ for  $i=1,3$, it follows that
$\det((\a\otimes \phi)(B_2))\ne 0$. It follows from Theorem
\ref{thm:Tu22} that $C_*\otimes_{\Z[\pi]}\F(t)^k$ is acyclic and
\[ \tau(M,\a,\phi)=\prod\limits_{i=1}^3 \det((\a\otimes \phi)(B_i))^{(-1)^{i+1}}.\]
Therefore $\tau(M,\a,\phi)=\phi(\tau(M,\a,\psi))$.

For the part (2), using similar arguments as above one can easily
show that if $\Delta_{M,\phi}^{\a,1} \ne 0$ then
$\phi\left(\Delta_{M,\psi}^{\a,1}\right) \ne 0$.

\end{proof}

\begin{proof}[{\bf Proof of Theorem \ref{mainthmalex}}]
Clearly Theorem \ref{mainthmalex} follows from Theorem
\ref{thm:turaev} and Lemmas \ref{lem:delta03}, \ref{lem:delta2}
(applied to $\psi:\pi_1(M)\to FH_1(M)$ and $\phi:\pi_1(M)\to \Z$)
and from Lemma \ref{lem:delta2multi}.
\end{proof}

\end{document}